\documentclass[12pt]{article}
\usepackage[nottoc]{tocbibind}
\usepackage[latin1]{inputenc}
\usepackage{fancyhdr}
\usepackage{verbatim}
\usepackage{indentfirst}
\usepackage{graphicx}
\usepackage{epstopdf}
\usepackage{setspace}
\usepackage{color}
\usepackage{newlfont}
\usepackage{amssymb}
\usepackage[intlimits]{amsmath}
\usepackage{latexsym}
\usepackage{amsthm}
\usepackage{amscd}
\usepackage{mathrsfs}
\usepackage{leftidx}
\usepackage[dvips,top=60pt,bottom=60pt,left=65pt,right=65pt]{geometry}
\usepackage{mathtools, array, multirow}
\usepackage{hyperref}
\usepackage{rotating}
\usepackage{mathtools}
\usepackage{dsfont}
\usepackage{import}
\usepackage{tikz, tikz-cd}
\usepackage{pgfplots}
\usepgfplotslibrary{fillbetween}
\pgfplotsset{compat=1.15}
\usepackage{thm-restate}

\numberwithin{equation}{section}

\theoremstyle{plain}                    
\newtheorem{theorem}{Theorem}[section]
\newtheorem*{theorem*}{Theorem}
\newtheorem{lemma}[theorem]{Lemma}    
\newtheorem{proposition}[theorem]{Proposition}
\newtheorem{corollary}[theorem]{Corollary}
\theoremstyle{definition}

\newtheorem{remark}[theorem]{Remark} 

\newcommand{\todo}[1]{\textcolor{red}{\bf \boxed{TODO} #1}\PackageWarning{TODO:}{#1!}}
\newcommand{\comm}[1]{\textcolor{blue}{\bf \boxed{COMMENT} #1}\PackageWarning{TODO:}{#1!}}
\newcommand{\ques}[1]{\textcolor{blue}{\bf \boxed{QUESTION} #1}\PackageWarning{TODO:}{#1!}}

\newcommand{\Z}{\mathbb{Z}}

\newcommand{\R}{\mathbb{R}}
\newcommand{\C}{\mathbb{C}}
\newcommand{\Hh}{\mathbb{H}}

\newcommand{\calF}{\mathcal{F}}

\newcommand{\vecv}{\text{\boldmath$v$}}

\newcommand{\lpar}{\left(}
\newcommand{\rpar}{\right)}
\newcommand{\lbr}{\left\{}
\newcommand{\rbr}{\right\}}
\newcommand{\labs}{\left|}
\newcommand{\rabs}{\right|}

\newcommand{\SL}{\mathrm{SL}}
\newcommand{\sltr}{\mathrm{SL}_2 (\R)}
\newcommand{\sltz}{\mathrm{SL}_2 (\Z)}

\newcommand{\de}{\mathrm{d}}
\newcommand{\ind}{\mathds{1}}
\newcommand{\inv}{^{-1}}
\newcommand{\vol}{\mathrm{vol}}
\newcommand{\matr}[4]{\left( \begin{matrix} #1 & #2 \\ #3 & #4 \end{matrix} \right) }
\newcommand{\smatr}[4]{\bigr( \begin{smallmatrix} #1 & #2 \\ #3 & #4 \end{smallmatrix} \bigr) }

\newcommand{\bs}{\backslash}
\newcommand{\eps}{\varepsilon}

\newcommand{\tf}{\tilde{f}}
\newcommand{\im}{\mathrm{Im}}

\newcommand{\veps}{\varepsilon}
\newcommand{\psp}{\mathcal{R}}
\newcommand{\mainterm}{\mathcal M}
\newcommand{\err}{\mathcal E}
\newcommand{\numberthis}{\addtocounter{equation}{1}\tag{\theequation}}

\title{An Effective Slope Gap Distribution for Lattice Surfaces}
\author{Tariq Osman        \and
        Josh Southerland                       \and
        Jane Wang
}
\date{}
\begin{document}

\maketitle

\begin{abstract}
We prove an effective slope gap distribution result first for the square torus and then for general lattice translation surfaces. As a corollary, we obtain a dynamical proof for an effective gap distribution result for the Farey fractions. As an intermediate step, we prove an effective equidistribution result for the intersection points of long horocycles with a particular transversal of the horocycle flow in $\sltr/\Gamma$ where $\Gamma$ is a lattice. 
\end{abstract}


\section{Introduction}


In dynamical systems, especially ergodic theory, many theorems are qualitative in the sense that they identify when iterates of a dynamical system converge, but fail to quantify the rate of this convergence. In the following paper, we resolve an issue of this type: we provide an effective limit theorem for the gap distribution of the saddle connections on the square torus, which relates to the spacing statistics of Farey fractions, and then generalize the result to lattice translation surfaces. Before introducing the main results, we will describe the main objects studied, and explain what we mean by the gap distribution of saddle connections. 


A \textbf{translation surface} can be defined geometrically as a finite collection of disjoint polygons in the complex plane $\mathbb{C}$ with sides identified in parallel opposite pairs by translation. On a translation surface, there is a natural notion of length, area, and direction, inherited from the geometry of $\mathbb{C}$ as identified with $\mathbb{R}^2$. A translation surface can alternatively thought of as a pair $(X,\omega)$ where $X$ is a Riemann surface and $\omega$ is a holomorphic one-form. To go from the first definition to the second, we notice that our polygons can inherit the complex structure from $\mathbb{C}$ as well as the one-form $dz$. To go from the second definition to the first, one can show that every translation surface $(X,\omega)$ can by triangulated (see \cite{Wright} for more details). 

The geometry on a translation surface locally looks like Euclidean space $\mathbb{R}^2$ except at finitely many \textbf{cone points} that have cone angle $2 \pi n$ for some $n \geq 2$. There cone points correspond to the zeros of the one-form $\omega$. A \textbf{saddle connection} is a straight line segment connecting two cone points that does not pass through any other cone points. If $\gamma$ is an oriented saddle connection, then its \textbf{holonomy vector} is $\mathbf{v}_\gamma = \int_\gamma \omega$, the vector that describes the length and direction of $\gamma$. 

The set of saddle connection vectors 
\begin{equation} 
\label{eq:Lambda}
\Lambda := \{\mathbf{v}_\gamma : \gamma \text{ is a saddle connection of } (X,\omega)\}
\end{equation} 
is a discrete subset of $\mathbb{R}^2.$ One can then ask about the growth and randomness of the angles of $\Lambda$. 

Results of Masur (\cite{Masur88} and \cite{Masur90}) show that for almost every translation surface, the set $\Lambda \cap B(0,R)$ of saddle connections of length at most $R$ has quadratic upper and lower bounds. By results of Vorobets (\cite{Vorobets}), it is also known that for almost every $(X,\omega)$ with respect to a natural measure on the space of translation surfaces called the Masur-Veech volume (see, for example, \cite{Zorich} for more details), the angles of the saddle connections equidistribute in the circle. That is, as $R \rightarrow \infty$, the proportion of saddle connections with length $\leq R$ with angle in an interval $I \subset S^1$ converges to the length of $I$. Here, we normalize $S^1$ to have total length $1$. 

We note that equidistribution of a given sequence on $S^1$ is akin to the sequence being generated by a sequence of independent, uniform random variables on $S^1$. To gain a deeper understanding as to just how random a given sequence is, we may study the spacing statistics of the sequence. 

For dynamical reasons, it is often easier to consider the spacing statistics of slopes of saddle connections, rather than the angles. Following a convention of Athreya, Chaika, and Leli\'{e}vre (\cite{Athreya-Chaika-Lelievre}), we let 
\begin{equation} 
\label{eq:Lambda_R}
\Lambda_R := \{v=(a,b)\in \Lambda: 0 \leq b \leq a \leq R\}.
\end{equation} 
That is, $\Lambda_R$ denotes the vectors in $\Lambda$ with a positive $\leq R$ horizontal component and slope between $0$ and $1$. We will be concerned with the set of slopes of this set as $R \rightarrow \infty$.

To assess the randomness of saddle connection directions on $(X,\omega)$, we wish to understand the \textbf{limiting slope gap distribution}. To define the slope gap distribution of $(X,\omega)$, we first order the $N(R)$ unique slopes of $\Lambda_R$ in increasing order: 

$$0 \leq s_R^{(0)} < s_R^{(1)} < \ldots < s_R^{N(R)-1} \leq 1.$$

Here $N(R)$ grows quadratically. We then consider the set of \textbf{renormalized gaps} 

\begin{equation}
    \label{eq:G_R}
\mathbb{G}_R := \{R^2( s_R^{(i)}-s_R^{(i-1)}) : 1 \leq i \leq N(R)\},
\end{equation} 
where the renormalization by $R^2$ is to counter that the expected gap size decays like $R^2$. 

We say that $F:[0,\infty) \rightarrow [0,\infty)$ is the \textbf{slope gap distribution} of $(X,\omega)$ if, for every $a < b$, we have that $$\lim_{R\rightarrow \infty} \frac{|\mathbb{G}_R \cap (a,b)|}{N(R)} = \int_a^b F(x) \, dx.$$ 

It is known by work of Athreya and Chaika (\cite{Atheya-Chaika}) that for almost every translation surface $(X,\omega)$ with respect to the Masur-Veech volume, the slope gap distribution of $(X,\omega)$ exists. If the sequence of slopes of saddle connections ordered by length were independent and uniformly distributed in $[0,1]$, then one would expect the slope gap distribution to be exponentially distributed. It is currently unknown if there exist any translation surfaces whose slope gap distribution is exponential, and the study of slope gap distributions of general translation surfaces is difficult. 

Athreya and Cheung (\cite{Athreya-Cheung-poincare-section}) computed the slope gap distribution of the square torus, which is also the gap distribution of the Farey fractions. Work of Uyanik and Work (\cite{Uyanik-Work}) as well as Kumanduri, Sanchez, and Wang (\cite{Kumanduri-Sanchez-Wang}) gives an algorithm to compute the slope gap distributions of an important subclass of translation surfaces with many symmetries called \textbf{lattice surfaces}, also known as \textbf{Veech surfaces}. This work will be described in more detail in Section \ref{sec:background}. 

While the existence of slope gap distributions is known in a few cases, to the best of our knowledge, there are no known results about how quickly renormalized slope gaps converge to the gap distribution. Our contribution is to make effective the convergence to the gap distribution in these two cases: the square torus and Farey fractions, and then general lattice surfaces. 

For the following two theorems, let $N(R)$ denote the number of unique slopes of saddle connections in $\Lambda_R$ (defined in Equation \ref{eq:Lambda_R}) and $\mathbb{G}_R$ denote the renormalized slope gaps of $\Lambda_R$, as defined earlier in the introduction.

\begin{restatable}[Effective gap distribution (torus and Farey fractions)]{theorem}{torusgaps}
\label{thm:torus-gaps}
    Let $F(x)$ be the slope gap distribution of the square torus with one marked point, which is also the limiting gap distribution of the Farey fractions. Let $\mathbb{G}_R$ be the renormalized gaps of the slopes $\Lambda_R$ for the square torus, or the renormalized gaps of the Farey fractions of denominator $\leq R$. Then, there exists a constant $C > 0$ such that for any $0 \leq a < b < \infty$, $$\left| \frac{|\mathbb{G}_R \cap (a,b)|}{N(R)} - \int_a^b F(x) \, dx \right| \leq C\log(R)R^{-\frac{1}{15}}$$.
\end{restatable}

\begin{remark} We note that renormalized gaps in the sets of Farey fractions $$\mathcal{F}(Q) := \{ \text{reduced fractions } \frac{p}{q} : (q,p) \in \mathbb{Z}^2, 0 < q \leq Q\}$$ are known to converge to a function called \textbf{Hall's distribution} (see \cite{Athreya-Cheung-poincare-section} and Section \ref{sec:torus-background} for more details). The Farey fractions $\mathcal{F}(Q)$ are also exactly the slopes $\Lambda_Q$ of the square torus. In this way, it can be seen that the slope gap distribution of the square torus is the gap distribution of the Farey fractions, and the rate of convergence in Theorem \ref{thm:torus-gaps} also applies to the convergence of the Farey fraction gaps to Hall's distribution. 
\end{remark}

\begin{restatable}[Effective gap distribution (lattice surface)]{theorem}{veechgaps}
\label{thm:Veech-gaps} 
    Let $F(x)$ be the slope gap distribution of the lattice surface $(X,\omega)$. Then, for any $0 \leq a < b < \infty$, $$\left| \frac{|\mathbb{G}_R \cap (a,b)|}{N(R)} - \int_a^b F(x) \, dx \right| \leq \begin{cases} 
      C  \log(R) R^{-\frac{1}{15}} & \text{if } \Gamma \text{ is tempered}\\
      C(s)  R^{-\frac{1}{15}(1-s)} & \text{if } \Gamma \text{ is non-tempered.}
\end{cases} $$ where the constant $C$ depends on the choice of lattice surface, and if $\Gamma$ is not tempered, $s$, where $\frac{1-s^2}{4}$ is the bottom of the spectrum of the hyperbolic Laplacian on $\Hh/\Gamma$.
\end{restatable}

For our purposes, a \textbf{tempered lattice} is one such that $\frac{1}{4}$ is the bottom of the spectrum of the hyperbolic Laplacian on $\Hh/\Gamma$, noninclusive of the $0$, whereas a non-tempered lattice is one such that the bottom of the spectrum of the hyperbolic Laplacian on $\Hh/\Gamma$ is strictly less than $\frac{1}{4}$. The spectral parameter $s$ in Theorem \ref{thm:Veech-gaps} quantifies the spectral gap.

\begin{remark} The proofs of Theorems \ref{thm:torus-gaps} and \ref{thm:Veech-gaps} rely on the effective equidistribution results of Theorems \ref{thm:eff-equidistribution-on-transversal} and \ref{thm:eff-equidistribution-on-transversal-Veech-surfaces}. As Remark \ref{rmk:exponents-for-effective-equidistribution} indicates, the rate of decay in these latter results is actually better than stated. Consequently, the exponents in Theorems \ref{thm:torus-gaps} and \ref{thm:Veech-gaps} are also better than stated. One can change the $\frac{1}{15}$ to $\frac{1}{3}$ in Theorem \ref{thm:torus-gaps} for the torus. In Theorem \ref{thm:Veech-gaps} for general lattice surfaces, when $(X,\omega)$ is periodic under the unstable horocycle flow, one can change $\frac{1}{15}$ to $\frac{1}{4n} - \varepsilon$ for any $\varepsilon > 0$, where $n$ is the smallest even integer greater than $\frac{1}{1-s}$ and $\frac{1 + s^2}{4}$ is the bottom of the spectrum of the hyperbolic Laplacian on $\Hh/\Gamma$. In the case where $(X,\omega)$ is not periodic under the unstable horocycle flow, one can change the $\frac{1}{15}$ to $\frac{1}{9}$ in Theorem \ref{thm:Veech-gaps}. 
\end{remark}


The proof of these theorems rests on the following effective equidistribution results for the intersection points of a family of long horocycles with a particular Poincar\'e section $\Omega$, which is a transversal to the horocycle flow. The horocycle flow, the transversal $\Omega$, their connection to gap distributions, and the natural measure $m$ on $\Omega$ will be discussed in Section \ref{sec:background}. Briefly, if $\Gamma \subset \sltr$ is the group of symmetries of $(X,\omega)$, then $\Omega \subset \sltr/\Gamma$ is the set of translation surfaces in the $\sltr$ orbit of $(X,\omega)$ that has a short (length $\leq 1$) horizontal saddle connection.


\begin{restatable}[Effective equidistribution on the section (torus)]{theorem}{effectivetorus}
\label{thm:eff-equidistribution-on-transversal}
Let $H_L$ be a the set of points corresponding to a segment of length $L$ of an (unstable) periodic horocycle trajectory. Let $H_{L(t)} = g_t^{-1}H_L$, where $g_t$ is the geodesic flow for time $t$, and $L(t)$ denotes the length of $H_{L(t)}$. Let $\rho_{L(t)} = \frac{1}{\lvert H_{L(t) \cap \Omega} \rvert}\sum_{x \in H_{L(t)} \cap \Omega} \delta_{x}$ and let $m$ be the measure on the section $\Omega$. For any compactly supported and bounded $f \in L^2(\Omega)$, 
\begin{equation*}
\left| \rho_{L(t)}(f) - m(f) \right| \leq  C \| f \|_2 \log(L(t)) L(t)^{-\frac{1}{30}}
\end{equation*}
where the constant $C$ depends on the supremum of $f$ and the geometry of the section.
\end{restatable}

For general lattice surfaces, we have a similar effective equidistribution results with worse bounds.

\begin{restatable}[Effective equidistribution on the section (lattice surface)]{theorem}{effectiveveech}
\label{thm:eff-equidistribution-on-transversal-Veech-surfaces} 
Let $f \in L^2(\Omega)$ be a compactly supported and bounded function. Then, with the same notation as defined in Theorem \ref{thm:eff-equidistribution-on-transversal},

\begin{equation*}
\left| \rho_{L(t)}(f) - m(f) \right| \leq \begin{cases} 
      C \| f \|_2 \log(L(t)) L(t)^{-\frac{1}{30}} & \text{if } \Gamma \text{ is tempered}\\
      C(s) \| f \|_2 (L(t))^{-\frac{1}{30}(1-s)} & \text{if } \Gamma \text{ is non-tempered.}
\end{cases}
\end{equation*}
where the constant $C$ depends on the supremum of $f$, the geometry of the section $\Omega$, and if $\Gamma$ is not tempered, $s$, where $\frac{1-s^2}{4}$ is the bottom of the spectrum of the hyperbolic Laplacian on $\Hh/\Gamma$.
\end{restatable}

\begin{remark} 
\label{rmk:exponents-for-effective-equidistribution}
The proofs of Theorems \ref{thm:eff-equidistribution-on-transversal} and \ref{thm:eff-equidistribution-on-transversal-Veech-surfaces} rely on the effective equidistribution of long horocycle segments. We have provided a dynamical proof of this well-known result in the appendix, and we use the rates that we compute in deducing our Theorems \ref{thm:eff-equidistribution-on-transversal} and \ref{thm:eff-equidistribution-on-transversal-Veech-surfaces}. However, using the best results in the literature, in the case of closed horocycles, one can change the $\frac{1}{30}$ to $\frac{1}{6}$ in Theorem \ref{thm:eff-equidistribution-on-transversal}. In Theorem \ref{thm:eff-equidistribution-on-transversal-Veech-surfaces}, for closed horocycles, one can change $\frac{1}{30}$ to $\frac{1}{8n} - \varepsilon$ for any $\varepsilon > 0$, where $n$ is the smallest even integer greater than $\frac{1}{1-s}$ and $\frac{1 + s^2}{4}$ is the bottom of the spectrum of the hyperbolic Laplacian on $\Hh/\Gamma$. Moreover, in the case of non-closed horocycles, one can change the $\frac{1}{30}$ to $\frac{1}{18}$ in both Theorem \ref{thm:eff-equidistribution-on-transversal} and \ref{thm:eff-equidistribution-on-transversal-Veech-surfaces}. The reader may follow Remarks \ref{rmk:long-horocycle-improvements}. \ref{rmk:rates-for-torus}, and \ref{rmk:rates-for-Veech-surfaces} to track these rates, and for additional information and references. 
\end{remark}

\begin{remark} The bounds for lattice surfaces are, in general, worse because for many lattice surfaces, $\Gamma$ is not a tempered lattice. Our arguments rely crucially the action of $\sltr$ on $\sltr/\Gamma$, and as a result, the appearance of a complementary series representation is an obstruction: the $SL_2(\R)$-action commutes with the Casimir operator on $\sltr/\Gamma$, and hence iterates of the action fall prey to the smallest element in the spectrum of the operator. It is well-known that this eigenvalue agrees with the bottom of the spectrum of the hyperbolic Laplacian on $\Hh/\Gamma$ since complementary series representations are spherical.
\end{remark}



\textbf{Related Literature.} In \cite{Athreya-Cheung-poincare-section}, Athreya and Cheung connected the gap distribution of the Farey fractions back to the slope gap distribution of the square torus. They pioneered a method that used renormalization techniques for the horocycle flow to compute the slope gap distribution. In \cite{Atheya-Chaika}, Athreya and Chaika analyzed the slope gap distributions of generic translation surfaces and also showed that lattice surfaces have \textit{no small gaps}. That is, the $\liminf$ of the renormalized gaps is bounded away from zero if and only if a surface is a lattice surface. More examples of slope gap distributions of lattice surfaces were computed by Athreya, Chaika, and Leli\`{e}vre \cite{Athreya-Chaika-Lelievre} for the golden L, by Uyanik and Work \cite{Uyanik-Work} for the regular octagon, and by Berman et al \cite{Berman-2n} for general $2n$-gons. In \cite{Uyanik-Work}, Uyanik and Work also proposed an algorithm for computing the slope gap distributions of arbitrary lattice surfaces. That algorithm was later improved upon by Kumanduri, Sanchez, and Wang in \cite{Kumanduri-Sanchez-Wang}. In \cite{Sanchez}, Sanchez proves properties about the gap distributions of double-slit tori, a family of non-lattice surfaces.

The proofs of Theorems \ref{thm:eff-equidistribution-on-transversal} and \ref{thm:eff-equidistribution-on-transversal-Veech-surfaces} rests on the effective equidistribution horocycle translates, under the action of the geodesic flow. In the case where the horocycle is closed, then the effective result is due to Sarnak \cite{Sarnak}. In \cite{Strombergsson}, the effective equidistribution of generic horocycle arcs are considered. Theorems \ref{thm:athreya-cheung} and \ref{thm:eff-equidistribution-on-transversal-Veech-surfaces} are derived from Theorem \ref{thm:eff-equid-long-hor-mixing}, which is a consequence of Theorem \ref{thm:decay-of-matrix-coeff}, and the ``Margulis Thickening Technique". We recount the proof of Theorem \ref{thm:eff-equid-long-hor-mixing} in Appendix \ref{appendix:proof-of-effective-equidist}, making explicit the dependence of the implied constant on the test function. We remark that the error term in Theorem \ref{thm:eff-equid-long-hor-mixing} is notably worse than those in the results of \cite{Sarnak} and \cite{Strombergsson}, but allows for a wider class of test functions. 

Theorems concerning the equidistribution of expanding translates of horospheres have been shown to hold in very general contexts, see for instance \cite{Edwards,Kleinbock-Margulis-09, LMW-Effective, strombergsson-effective-lift}, and references therein. Such theorems have found application in a wide variety of problems across modern mathematics, such as the Berry-Tabor conjecture on flat tori \cite{LMW-flat-torus},  the gap distribution of $\sqrt n$ modulo $1$ \cite{Browning-Vinogradov}, the asymptotic distribution of Frobenius numbers \cite{Marklof-frobenius, Li-Frobenius}, and Apollonian circle packings \cite{Athreya-Cobeli-Zaharescu}, to name a few.




\textbf{Structure of the proofs of Theorems \ref{thm:eff-equidistribution-on-transversal} and \ref{thm:eff-equidistribution-on-transversal-Veech-surfaces}.} Recall that $H_{L(t)}$ is a geodesic flow push of a horocycle segment $H_L$ in $\sltr/\Gamma$. The main technical step in proving an effective gaps result (Theorems \ref{thm:torus-gaps} and \ref{thm:Veech-gaps}) is to understand the rate at which the counting measure $\rho_{L(t)}$ of $H_{L(t)} \cap \Omega$ converges to Lebesgue measure $m$ on the transversal $\Omega \subset \sltr/\Gamma$ over a suitably regular class of test functions. More precisely, for sufficiently regular $f$, we estimate 
\begin{equation}
    \labs \rho_{L(t)}(f) - m(f)\rabs.\label{eq:rho-to-m}
\end{equation}

Given a function $f$ on $\Omega$, we define a particular `thickening' of $f$ to a function $\tf$, defined on $\sltr / \Gamma$. 
The function $\tf$ is chosen so that it is constant in unstable horocycle direction. As we will see in Section \ref{sec:gaps-algorithm}, $\Omega$ breaks into finitely many pieces $\Omega_i$ where there are natural local $(a,b)$ coordinates, and the thickening width, $w$, can be chosen depending on $(X,\omega)$ so that $w$ is less than the minimum horocycle return time on $\Omega$.

We then define
\begin{align}
\nu_{L(t)} &:= \text{length measure on } H_{L(t)} \subset \sltr/\Gamma, \text{ normalized so the total measure is }1 \\
\mu & := \text{Haar measure on } \sltr/\Gamma, \text{normalised to be a probability measure.}
\end{align}

By multiple applications of the triangle inequality, finding an upper bound for \eqref{eq:rho-to-m} amounts to finding upper bounds for $| \rho_{L(t)}(f) - \nu_{L(t)}(\tilde{f})|$, $| \nu_{L(t)}(\tilde{f}) - \mu(\tilde{f})|$, and $|\mu(\tilde{f}) - m(f) |$, respectively.

For $|\rho_{L(t)}(f) - m(f)|$, we use the relation between the length of $H_{L(t)}$ and $|H_{L(t)} \cap \Omega|$ as proven in Section \ref{sec:comparisons}. This relation is derived from results about counting lattice points in certain regions of the plane. An upper bound for $|\rho_{L(t)}(f) - \nu_{L(t)}(\tilde{f})|$ is obtained using the effective equidistribution of horocycle translates. In order to apply the appropriate theorems, we must first smooth $\tf$. This procedure requires various technical lemmas, carried out in Section \ref{sec:effective-transversal}. Finally, $|\mu(\tilde{f}) - m(f) |=0$ as shown in Proposition \ref{prop:mu-and-m}. We note that there is a tension between when the decay rates of $|\rho_{L(t)}(f) - m(f)|$ coming from lattice counting methods or the decay rates of $|\rho_{L(t)}(f) - \nu_{L(t)}(\tilde{f})|$ coming from the effective equidistribution of long horocycles contribute the dominant term in the error rate of our effective gaps theorems. This is commented on more in Section \ref{sec:effective-transversal}.

\textbf{Outline of the Paper.} 
In Section \ref{sec:background}, we summarize the relevant existing work on the relation between slope gap distributions and the horocycle flow and how this can be used to compute the slope gap distribution of the Farey fractions and the square torus. This section also contains the proofs of some key propositions as well as our main effective gaps results, Theorems \ref{thm:torus-gaps} and \ref{thm:Veech-gaps}, assuming our effective equidistribution results, Theorems \ref{thm:eff-equidistribution-on-transversal} and \ref{thm:eff-equidistribution-on-transversal-Veech-surfaces}. The remaining sections are devoted to proving these effective equidistribution results. In Section \ref{sec:comparisons}, we prove a key relationship between the lengths of long horocycle segments and the number of times that they intersect the transversal $\Omega$. We prove these bounds first for the square torus and then for general lattice surfaces. Section \ref{sec:effective-transversal} then contains the proof Theorems \ref{thm:eff-equidistribution-on-transversal} and \ref{thm:eff-equidistribution-on-transversal-Veech-surfaces}, which are effective equidistribution results on the transversal $\Omega$ to the horocycle flow. 

We note that there are two main goals of this paper: to prove an effective gaps result for the square torus and the Farey fractions, and to prove a more general effective gaps result for all Veech translation surfaces. A reader interested primarily in Farey fraction gaps could skip Sections \ref{sec:gaps-algorithm}, \ref{sec:comparisons-veech}, and \ref{sec:veech-transversal-equidistribution},  whereas a reader interested also in the effective gaps of lattice surfaces would read the whole paper.

\textbf{Acknowledgements.} We would like to thank Anthony Sanchez for suggesting this problem to us, and for interesting initial discussions. This project began from conversations held at the Summer School on Renormalization and Visualization for Packing, Billiards, and Surfaces at CIRM, Centre Internationale de Recontres Math\'ematiques. We would like to thank CIRM for hosting the conference and give a special thank you to the Chaire Jean Morlet program, the Jean Morlet Chair being Jayadev Athreya at the time of the confernce. We would also like to thank Jayadev Athreya for funding the second and third author's trips to CIRM for the conference with NSF grant DMS 2333366. 
\section{Background and Effective Gaps}
\label{sec:background}

In this section, we introduce the necessary background on slope gap distributions and prove our effective gap distribution theorems. Most of the content in Sections \ref{sec:horocycle-background} through \ref{sec:gaps-algorithm} is not new, excepts that in Section \ref{sec:UW-proof}, we note that a key gap distribution proof from Uyanik and Work (\cite{Uyanik-Work}) does not work as stated for all lattice surfaces. We comment on why and remedy this. Section \ref{sec:return-time-function-bounds} contains two new propositions needed to prove our main theorems, and Section \ref{sec:effective-gap-proofs} contains the proofs of our main effective gaps theorems. 

\subsection{From the horocycle flow to slope gap distributions }
\label{sec:horocycle-background}

There is a natural action of $\sltr$ on the space of translation surface coming from the action of $\sltr$ on $\mathbb{R}^2$. If $M \in \sltr$ and $(X,\omega)$ is a translation surface made by gluing the polygons $P_1, \ldots, P_k$, then $M \cdot (X,\omega)$ is the translation surface made by gluing the polygons $MP_1, \ldots, MP_k$ in the same gluing pattern. Here, $M$ acts on each $P_i \subset \mathbb{R}^2$ via the linear action of $M$ on the plane. 

Sometimes, $M(X,\omega)$ and $(X,\omega)$ are equivalent after cut and paste operations. In this case, $M$ is said to be in the \textbf{Veech group} $\Gamma = \SL(X,\omega) \subset \sltr$. That is, the Veech group is the subgroup of all matrices that stabilize $(X,\omega)$. When $\SL(X,\omega) \subset \sltr$ is a lattice (that is, when $\sltr/\SL(X,\omega)$ has finite volume), $(X,\omega)$ is said to be a \textbf{lattice surface}. Lattice surfaces are a measure zero set in the space of all translation surfaces, but are important examples of translation surfaces, in part because they have nice dynamical properties. In particular, it is known that the $\sltr$ orbit of a translation surface $(X,\omega)$ is closed in its \textbf{stratum} of the moduli space, which includes all translation surfaces with the same cone angle data, if and only if $(X,\omega)$ is a lattice surface. This closed orbit property makes computing the slope gap distributions of lattice surfaces more tractable than for general translation surfaces, where the orbit closure is higher dimensional and dynamics on the orbit closure is more complicated. For the rest of this paper, all translation surfaces $(X,\omega)$ are assumed to be lattice surfaces. 

Athreya and Chaika (\cite{Atheya-Chaika}) have shown that the slope gap distributions of lattice surfaces have no support near zero and therefore cannot be exponential. However, we now know how to compute the slope gap distribution. An algorithm to compute this limiting gap distribution for lattice surfaces was first proved by Uyanik and Work (\cite{Uyanik-Work}), and then was later improved upon by Kumanduri, Sanchez, and Wang (\cite{Kumanduri-Sanchez-Wang}). This algorithm relies on renormalization techniques that allow us to relate the slope gap distribution of a lattice surface $(X,\omega)$ to the horocycle flow on the $\sltr$ orbit of $(X,\omega)$. These techniques were first exploited by Athreya and Cheung on the torus\cite{Athreya-Cheung-poincare-section}.


Let $h_s : = \begin{bmatrix} 1 & 0 \\ -s & 1 \end{bmatrix}$ denote the (unstable) horocycle flow and let $g_t := \begin{bmatrix} e^{t/2} & 0 \\ 0 & e^{-t/2}\end{bmatrix}$ denote the geodesic flow. Note that $h_s$ acts on translation surfaces via a vertical shear. Recall that $\Lambda$ is the set of holonomy vectors of $(X,\omega)$ Then, the space 
\begin{equation}
\label{eq:transversal}
\Omega := \{g\Gamma \in \sltr/\Gamma : g\Lambda \text{ contains a short (length $\leq 1$) horizontal vector}\}
\end{equation}
is a transversal for the horocycle flow, and $\sltr/\Gamma$ is a suspension over $\Omega$ where the height above a point is the return time of the $h_s$ flow to $\Omega$. Another way to think of $\Omega$ is as the set of surfaces $g(X,\omega)$ in the $\sltr$ orbit of $(X,\omega)$ that have a short horizontal saddle connection. $\Omega$ naturally inherits a measure $m$ from the Haar measure $\mu$ on $\sltr$, such that $dm \, ds = d \mu$, where $s$ is the horocycle direction. We can then define the return time function $R : \Omega \rightarrow \mathbb{R}_+$ by letting $R(x)$ be the first return time of $x \in \Omega$ back to $\Omega$ under the horocycle flow $h_s$. We also can define $T : \Omega \rightarrow \Omega$ by letting $T(x) = h_{R(x)}(x)$ to be the first return map of $x$ to $\Omega$ under the horocycle flow. We note that it also makes sense to extend $T$ to a map $\sltr/\Gamma \rightarrow \Omega$ by defining $T$ as $h_s(x)$ for the first $s>0$ for which $h_s(x) \in\Omega$.  

The following theorem is implicit in the paper of Uyanik and Work \cite{Uyanik-Work}. A proof sketch of this theorem will be given in Section \ref{sec:UW-proof}.

\begin{theorem}[\cite{Uyanik-Work}]
\label{thm:gaps} 
Let $(X,\omega)$ be a lattice surface and let $F: \mathbb{R}_+ \rightarrow \mathbb{R}_+$ be the distribution of the $h_s$ return time function. That is, $\int_a^b F(x) \, dx = m \{x \in \Omega : a \leq R(x) \leq b\}$. Then, $F$ is the slope gap distribution of every surface in the $\sltr$ orbit of $(X,\omega)$. 
\end{theorem}

With this theorem, the computation of a slope gap distribution for a lattice surface $(X,\omega)$ reduces to finding good coordinates for $\Omega$, understanding the return time function $R$ of the horocycle flow on $\Omega$ in these coordinates, and then finding the return time function $f$ using calculus. This process for general lattice surfaces was proposed in \cite{Uyanik-Work} and improved in \cite{Kumanduri-Sanchez-Wang}. A summary of this process is given in Section \ref{sec:gaps-algorithm}. 

\subsection{The square torus and Farey fractions}
\label{sec:torus-background}

An important, independently interesting, and illustrative example is that of the square torus,  created by identifying the opposite sides of a $1\times 1$ square and marking a single point. In this section, we highlight some details of the analysis of the square torus conducted in \cite{Athreya-Cheung-poincare-section}. The case of the square torus gives the gap distribution of the Farey fractions and introduces ideas that will appear again in the discussion of general lattice surfaces. 

Let $(X,\omega)$ be the square torus with one marked point and let 
\begin{equation}
\label{eq:Lambda_Q}
\Lambda_Q= \{\text{holonomy vectors of } (X,\omega) \text{ with slope $\frac{p}{q}$ in } [0,1] \text{ and $x$ coordinate in } (0,Q]\}.
\end{equation}
The slopes of the vectors in $\Lambda_Q$ are then exactly equal to the \textbf{Farey sequence} of order $Q$, sometimes also known as the \textbf{Farey fractions}:  

\begin{equation} 
\label{eq:Farey-seq}
\calF(Q) := \lbr \tfrac{p}{q} \in [0,1] : (q, p) \in \Z^2, \: 0 < q \leq Q\rbr. 
\end{equation}
Thus, letting $Q \rightarrow \infty$, the slope gap distribution of the square torus is exactly the limiting gap distribution of the Farey fractions. We note that with the ordering of coordinates $(q,p)$, the sets $\mathcal{F}(Q)$ and $\Lambda_Q$ from equations \ref{eq:Farey-seq} and \ref{eq:Lambda_Q} agree. 

In this case, it can be seen that the Veech group of $(X,\omega)$ is $\Gamma = \sltz$, and the moduli space of flat tori is given by $X_2 := \sltr/\sltz$.

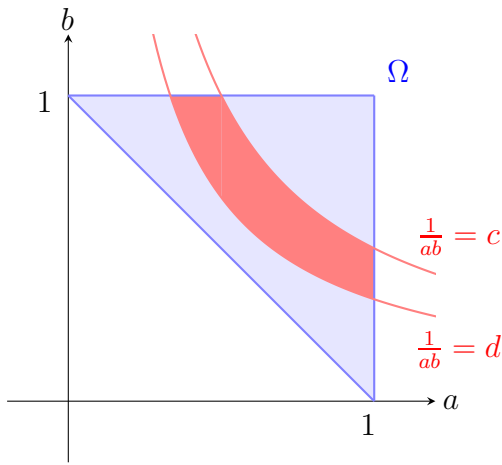
\begin{figure}
\centering
\begin{tikzpicture}
\begin{axis}[xmin=-0.2, xmax=1.2,
   		ymin=-0.2, ymax=1.2, axis lines = center, ticks = none, unit vector ratio=1 1 1]
    \addplot [domain=0:1, samples=100, name path=f, thick, color=blue!50]
        {1};
    
    \addplot [domain=0:1, samples=100, name path=g, thick, color=blue!50]
        {1-x};
        
    \addplot +[mark=none,color=blue!50,thick] coordinates {(1, 1) (1,0)};
    
    \addplot[blue!10, opacity=0.6] fill between[of=g and f, soft clip={domain=0:1}];
        
    \addplot [domain=0:1.25, samples=100, name path=j, thick, color=red!50]
        {1/(3*x)};

     \addplot [domain=0:1.25, samples=100, name path=k, thick, color=red!50]
        {1/(2*x)};

    \addplot[red!50, opacity=0.8] fill between[of=j and k, soft clip={domain=0.5:1}];

    \addplot[red!50, opacity=0.8] fill between[of=j and f, soft clip={domain=0.333:0.5}];
    
\end{axis}
\node[blue] at (5.2,5.2) {$\Omega$};
\node at (0.5,4.8) {1};
\node at (4.8,0.5) {1};
\node at (.8,5.9) {$b$};
\node at (5.9,.8) {$a$};
\node[red] at (6,3) {$\frac{1}{ab} = c$};
\node[red] at (6,1.5) {$\frac{1}{ab} = d$};
\end{tikzpicture}
\caption{The triangular transversal $\Omega$ in $(a,b)$ coordinates, along with the shaded region whose measure gives $\int_c^d F(x) \, dx$.}
\label{fig:Omega}
\end{figure}

Let $p_{a,b} := \begin{bmatrix} a & b \\ 0 & a^{-1} \end{bmatrix}.$ It was shown in \cite{Athreya-Cheung-poincare-section} that $\Omega$, as defined in Equation \ref{eq:transversal}, can be written explicitly in coordinates as

\begin{equation}
    \Omega := \{p_{a,b} \sltz : a, b \in (0,1], \: a + b > 1\} \subset X_2.\label{eq:poincare-section}
\end{equation}
This triangular transversal is depicted in Figure \ref{fig:Omega}. We note that when $\Omega$ is viewed as a subset of $T^1(\Hh)$, $\Omega = \{(z, \vecv) \in T^1(\Hh) : \im z \geq 1, \vecv = i\},$ (parameterizing $T^1(\Hh) = \Hh \times \{\vecv \in \C : \|\vecv\| = 1\}$).

Athreya and Cheung make explicit the return time map and first return map in the following theorem.

\begin{theorem}[\cite{Athreya-Cheung-poincare-section} Theorem 1.1 therein]
\label{thm:athreya-cheung}
    The set $\Omega$ as defined in \eqref{eq:poincare-section} is a Poincar\'e section (tranvsersal) for the action of the horocycle flow $h_s$ on $X_2$. Furthermore, the first return time $R(a,b)$ of the point $(a,b)$ in the coordinates for $\Omega$ is given by
    \begin{align}
        R (a, b) =  \frac{1}{ab}.
    \end{align}
    The first return map $T : \Omega  \to \Omega$ can be explicitly described by the BCZ map
    \begin{align}
        T(a, b) = \lpar b, -a + \left\lfloor \frac{1 + a}{b}\right\rfloor b\rpar.
    \end{align}
\end{theorem}

Then, by applying Theorem \ref{thm:gaps}, the gap distribution $f(x)$ of the square torus and of the Farey fractions satisfies that $$\int_c^d F(x) \, dx = m \left\{(a,b): a,b \in (0,1], a + b > 1, c < \frac{1}{ab} < d\right\}.$$ This latter region is the intersection of the region between two hyperbolas intersected with $\Omega$ (see Figure \ref{fig:Omega}).

Using calculus, one can explicitly compute the cumulative distribution function $\int_0^x F(t) \, dt$ and then differentiate to get the slope gap distribution $F(x)$. In this case, we find that $F(x)$ is \textbf{Hall's distribution}, a piecewise real-analytic function with two points of non-analyticity. More details can be found in \cite{Athreya-Cheung-poincare-section}.

\subsection{Relating gap distributions to the return time function}
\label{sec:UW-proof}



In this section, we sketch a proof of Theorem \ref{thm:gaps}, which relates the slope gap distribution to the return time function of the horocycle flow to the transversal $\Omega$, as defined in Equation \ref{eq:transversal}. We note that the proof of this result in Uyanik and Work's paper (\cite{Uyanik-Work}) contains a gap as written and only works as written for the case when $(X,\omega)$ has a vertical saddle connection. In our proof sketch, we highlight and fix this gap. The key dynamical result needed to prove this theorem is an equidistribution result about a family of long horocycles that are $g_t$ pushes of a horocycle segment. In latter sections, the effective version of this statement will be the key input in our proofs.

Recall that the $\sltr$ orbit of $(X,\omega)$ can be identified with $\sltr/\Gamma$. We let $H_L$ denote a horocycle segment of length $L$ in $\sltr/\Gamma$, based at some $g\Gamma$, usually corresponding with $(X,\omega)$. 
\begin{equation}
    \label{eq:gamma}
H_L := \{h_s(g\Gamma) \in \sltr/\Gamma: 0 \leq s \leq L\}.
\end{equation}
This $L$ is often chosen to be the period of $H_L$ if the horocycle is periodic, or $1$ if it is not periodic. 

We also define a $g_t$ push of this horocycle segment as 

\begin{equation}
    \label{eq:gt-push}
H_{L(t)} := g_t^{-1}H_L = \{g_t^{-1} h_s (g\Gamma) : 0 \leq s \leq L\} = \{h_sg_{t}^{-1} (g\Gamma) : 0 \leq s \leq L \cdot e^t\}.
\end{equation}

\begin{remark}
\label{rem:L(t)}
The latter equality above is because $g_t^{-1} h_s = h_{se^t} g_t^{-1}$ and shows that the length $L(t)$ of $g_t^{-1} H_L = H_{L(t)}$ is $e^tL$. 
\end{remark}

These $H_L(t)$ give a family of closed horocycle orbits or horocycle segments, depending on whether $H_L$ was closed or not respectively. We note that as $t \rightarrow \infty$, the length of $H_{L(t)}$ also goes to infinity.

\begin{proposition}
\label{prop:renormalized-equidistribution} 
Let $(X,\omega)$ be a lattice surface and $L > 0$. Then, as $t \rightarrow \infty$, the long horocycle segments $H_{L(t)}$ as defined in Equation \ref{eq:gt-push} equidistribute in the $\sltr$ orbit of $(X,\omega)$, which is identified with $\sltr/\Gamma$. 
\end{proposition}

\begin{proof}

Consider a lattice surface $(X,\omega)$, with Veech group $\Gamma = \SL(X,\omega)$. We can then identify the $\sltr$ orbit of $(X,\omega)$ up to equivalence with $\sltr/\Gamma$  by identifying $(X,\omega)$ with the coset $\Gamma \in \sltr/\Gamma$. Theorem \ref{thm:eff-equid-long-hor-mixing} is then an effective equidistribution statement for the segments $g_t^{-1}H_L = H_{L(t)}$. 
\end{proof}

With this proposition, we can prove the slope gap distribution theorem of Uyanik and Work. 

\begin{proof}[Proof of Theorem \ref{thm:gaps}.] 

Let $(X,\omega)$ be a lattice surface. Let $\Lambda_R$ be the set of saddle connections of slope between $0$ and $1$ for which the $x$ coordinate satisfies $0 < x \leq R$. To understand the renormalized slope gaps $\mathbb{G}_R$ of the set $\Lambda_R$, we let $t=2\log(R)$ and renormalize $(X,\omega)$ by $g_{t}^{-1} = \begin{bmatrix} 1/R & 0 \\ 0 & R \end{bmatrix}$ so that the saddle connection vectors of $g_{t}^{-1} (X,\omega)$ with a short (length $\leq 1$) positive horizontal component and slope between $0$ and $R^2$ are exactly $g_{t}^{-}\Lambda_R$. 

Then let $H_L$ be a length $L=1$ horocycle beginning at $(X,\omega)$ and $H_{L(t)} = g_t^{-1}H_L$ be a geodesic push of $H_L$. Since the horocycle flow preserves slope gaps, $\mathbb{G}_R$ is exactly the horocycle return times of the set of points $H_{L(t)} \cap \Omega$ to $\Omega$. 

We let $\rho_{L(t)}$ be the counting measure of $H_{L(t)} \cap \Omega$ on $\Omega$ and $m$ be the Lebesgue measure on $\Omega$.  By Proposition \ref{prop:renormalized-equidistribution}, $H_{L(t)}$  equidistributes in $\sltr/\Gamma$ as $t \rightarrow \infty$. Since $\Omega$ is a transversal to the $h_s$ flow and the Haar measure on $\sltr/\Gamma$ can be decomposed as $da \, db \, ds$ times a normalizing constant where $da \, db$ times a normalizing constant is Lebesgue measure in local $(a,b)$ coordinates on $\Omega$ and $ds$ is the length measure in the $h_s$ direction, $\rho_{L(t)} \rightarrow m$ as measures on $\Omega$ (for more details, see \cite{Athreya-Cheung-poincare-section}). 

Then, for for any $0 \leq a < b < \infty$, let $f_{a,b}: \Omega \rightarrow \mathbb{R}_+$ be the characteristic function on $R^{-1}(a,b)$, where $R : \Omega \rightarrow \mathbb{R}_+$ is the horocycle return time function to $\Omega$. It follows then that $\rho_{L(t)}(f_{a,b}) = \frac{|\mathbb{G}_R \cap (a,b)|}{N(R)}$ and $m(f_{a,b}) = \int_a^b F(x) \, dx = m\{x \in \Omega : a \leq R(x) \leq b\}$. Since $\rho_{L(t)} \rightarrow m$ as $t = 2 \log(R) \rightarrow \infty$, we have that  $\frac{|\mathbb{G}_R \cap (a,b)|}{N(R)} \rightarrow \int_a^b F(x) \, dx$, showing that $F(x)$ is the limiting gap slope gap distribution of $(X,\omega)$. 

\end{proof}

\begin{remark} We note that the reason that Uyanik and Work's proof implicitly relied on $(X,\omega)$ being vertically periodic is because they cited the equidistribution of long periodic horocycles (\cite{Sarnak}) to show that the $h_s$ intersection points with $\Omega$ equidistributed. For the $h_s$ orbit of $g_{-2\log(R)}(X,\omega)$ for $0 \leq s \leq R^2$ to be part of a periodic horocycle, $g_{-2\log(R)}(X,\omega)$ must be periodic under $h_s$. But in the case when $(X,\omega)$ is a lattice surface this occurs if and only if $g_{-2\log(R)}(X,\omega)$ and therefore $(X,\omega)$ is vertically periodic, which is also equivalent to $(X,\omega)$ having a vertical saddle connection (see, for example, \cite{Hubert-Schmidt}). 
\end{remark}

\begin{remark} The convention of Athreya-Chaika-Leli\'{e}vre was to let $\Lambda_R$ consist of saddle connections of slope between $0$ and $1$ for which the $x$ coordinate was positive and $\leq R$. The above proof sketch would also work if we defined $\Lambda_R$ to be to be those saddle connections of slope between $0$ and $C$ for any $C > 0$. In particular, in the periodic case, a natural choice would be to let $C$ be the period of $(X,\omega)$ under the horocycle flow. We note that that the slope gap distribution $f$ from Theorem \ref{thm:gaps} is independent of the choice of $C$. 
\end{remark}

\begin{remark} We note that this argument also shows that every element in the $\sltr$-orbit of $(X,\omega)$ also has the same limiting slope gap distribution, since the $h_s$ return time function is an invariant of the $\sltr$-orbit of $(X,\omega)$. 
\end{remark}

Theorem \ref{thm:gaps} tells us that to compute the gap distribution of a lattice surface, it suffices to understand the return time function of the horocycle flow $h_s$ to the transversal $\Omega$. A general framework for doing so is given in \cite{Uyanik-Work} and \cite{Kumanduri-Sanchez-Wang}. As stated in the introduction, the goal of this paper is to build on this work and to make effective the convergence of slope gaps to the limiting gap distribution, first for the square torus and then for general lattice surfaces. 

\subsection{The slope gap distribution algorithm}
\label{sec:gaps-algorithm}


In Theorem \ref{thm:gaps}, it was established that the slope gap distribution of a Veech translation surface $(X,\omega)$ can be computed if we understand the transversal $\Omega$ to the horocycle flow of elements in the $\sltr$ orbit of $(X,\omega)$ with a short horizontal saddle connection, as well as the horocycle return time function $R : \Omega \rightarrow \mathbb{R}_+$. In this section, we will outline how to parametrize $\Omega$ and some properties of the return time function $R$, and provide some related intuition. We note that many of the details in this section are not necessary for understanding the proofs of our theorems, but we provide them for completeness (for even more detail, see \cite{Uyanik-Work} and \cite{Kumanduri-Sanchez-Wang}). The main takeaway from this section is that the transversal $\Omega$ has natural coordinates in which it breaks up into finitely many polygonal pieces on which the return time function is real analytic. This will allow us to prove a return time bound in Proposition \ref{prop:horocycle-return-time-bounds} that will be necessary in the proofs of the main theorems.

We suppose that $(X,\omega)$ is a Veech translation surface with $n < \infty$ cusps. That is, the Veech group $\Gamma = \SL(X,\omega)$ has $n$ conjugacy classes of maximal parabolic subgroups, with each subgroup contributing a cusp to the quotient $\sltr/\Gamma$. We let $\Gamma_1, \ldots, \Gamma_n$ be representatives of these conjugacy classes. 

If $\Lambda(X,\omega)$ is the set of all saddle connection vectors of $(X,\omega)$, then it can be shown that 
$$\Lambda(X,\omega) = \bigcup_{i=1}^m (\Gamma v_i)$$
is the disjoint union of a finite number of $\Gamma$ orbits $\Gamma v_i$ for vectors $v_i \in \mathbb{R}^2$ (see \cite{Athreya-Chaika-Lelievre}). Here, $m \geq n$, where $n$ is the number of cusps of $\sltr/\Gamma$, which is also equal to the number of conjugacy classes of  maximal parabolic subgroups of the Veech group $\Gamma$. For each such parabolic subgroup, there are one or more $v_i$ given by the saddle connections in the eigendirection of the infinite cyclic generator of $\Gamma_i$. 

As a result of this decomposition of the saddle connections $\Lambda(X,\omega)$, the transversal $\Omega$ of surfaces with a short horizontal saddle (as defined in Equation \ref{eq:transversal}) can be naturally broken up into $n$ disjoint pieces, one for each cusp. 

We have that $\Gamma_i \cong \mathbb{Z} \oplus \mathbb{Z}/2$ or $\mathbb{Z} $ depending on whether $-I$ is in or  not in $\Gamma = \SL(X,\omega)$ respectively. We can choose a generator $P_i$ of the infinite cyclic factor that has eigenvalue $\pm 1$, and let $v_i$ be the shortest saddle connection vector of $(X,\omega)$ that is an eigenvector of $P_i$. Then, there exists a $C_i \in \sltr$ such that $$C_iP_iC_i^{-1} = \begin{bmatrix}\pm 1 & \alpha_i \\ 0 & \pm 1 \end{bmatrix},$$ where the sign of $\pm1$ matches with the sign of the eigenvalue of $P_i$ and $C_i  v_i = \begin{bmatrix} 1 \\ 0 \end{bmatrix}.$

\begin{figure}[ht]
\centering
\begin{tikzpicture}
\draw [fill = red!25] (0,1) -- (3,-0.5) -- (3,-4) -- (0,1);
\node at (4.5,0) {$a$};
\node at (0,2.5) {$b$}; 
\node at (-.3, 1) {$\frac{1}{y_0}$};
\draw[thick, <->] (-1, 0) -- (4,0); 
\draw[thick, <->] (0,-4) -- (0,2); 
\draw (3,-.1) -- (3,.1); 
\node at (3,.3) {$1$}; 

\node at (2, .8) {$b = \frac{-x_0}{y_0} a + \frac{1}{y_0}$}; 

\node at (.2, -2.7) {$b = (\frac{-x_0}{y_0}-n) a + \frac{1}{y_0}$}; 
\end{tikzpicture}
\caption{A parametrization of $\Omega_i$, where $n = \alpha_i$ or $2\alpha_i$ depending on whether $-I \in \SL(X,\omega)$ and the eigenvalues of $P_i$.}
\label{fig:Omega_i}
\end{figure}
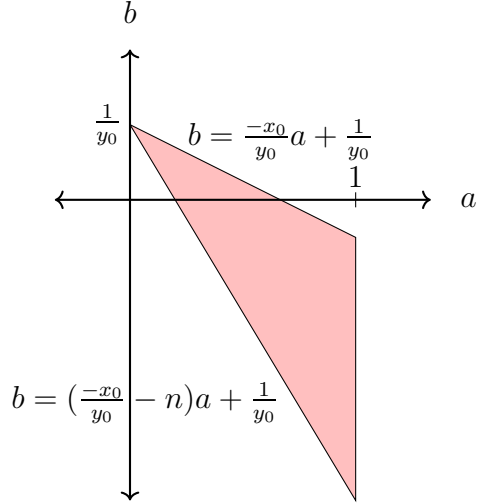

It follows that $C_i (X,\omega)$ is in the transversal $\Omega$ because $C_i v_i$ is a short horizontal saddle connection. The idea to parametrize the piece of $\Omega$ corresponding to $\Gamma_i$ is to notice that for a family of matrices $M_{a,b} \in \sltr$ (up to equivalence by $\Gamma$), $M_{a,b}\cdot C_i (X,\omega)$ has a short horizontal saddle connection $M_{a,b}\cdot C_i v_i$. After making some careful choices about which representative matrices $M_{a,b}$ of $\sltr/\Gamma$ to choose, it was shown in \cite{Kumanduri-Sanchez-Wang} that $\Omega$ could be parametrized by pieces $\Omega_i$ for each parabolic subgroup $\Gamma_i$ as follows: 

Consider $C_i(X,\omega)$ with $C_i$ as defined earlier in this section. $C_i(X,\omega)$ has a horizontal saddle connection $C_i v_i = \begin{bmatrix} 1 \\ 0 \end{bmatrix}$ and therefore breaks up into horizontal cylinders (see \cite{Hubert-Schmidt}). Let $y_0 > 0$ be the shortest cylinder height, and let $x_0 > 0$ be the shortest horizontal component of a saddle connection with height $y_0$. Then, $\Omega_i$ is defined in three cases: 
\begin{enumerate}
    \item When $-I \not \in \SL(X,\omega)$, $$\Omega_i := \{(a,b) \in \mathbb{R}^2: 0 < a \leq 1, \frac{1-x_0a}{y_0} - \alpha_i a \leq b \leq \frac{1-x_0a}{y_0}\}.$$
    \item When $-I \in \SL(X,\omega)$ and the eigenvalue of $P_i$ is $1$, then $\Omega_i$ consists of the triangle from case 1 as well as $-1$ times the triangle. 
    \item When $-I \in \SL(X,\omega)$ and the eigenvalue of $P_i$ is $-1$, $\Omega_i$ is the triangle from case 1 with $\alpha_i$ replaced by $2\alpha_i$. 
\end{enumerate}
The shape of these $\Omega_i$ pieces is depicted in Figure \ref{fig:Omega_i}. Then, $\Omega$ consists of a disjoint union of the pieces $$\Omega_i^M := \{M_{a,b} C_i (X,\omega) : (a,b) \in \Omega_i\}.$$ The $\Omega_i$ thus give computationally nice local coordinates for each piece of the transversal $\Omega$. We note that the vertex $(0,1/y_0)$ of each of these triangles is not included in $\Omega_i$ and a neighborhood of $(0,1/y_0)$ corresponds to the intersection of $\Omega_i$ with neighborhood of a cusp in $\sltr/\Gamma$. 

\begin{remark} The parametrization described here from Kumarduri-Sanchez-Wang (\cite{Kumanduri-Sanchez-Wang}) is modification of the parametrization from Uyanik-Work (\cite{Uyanik-Work}), chosen so that each piece $\Omega_i$ is guaranteed to break up into finitely many polygonal pieces, each with a piecewise real analytic return time function for the horocycle flow $h_s$. This generically did not happen for the parametrization chosen by Uyanik and Work. This finiteness will be useful when we prove properties of the return time function of $h_s$ to $\Omega$. 
\end{remark}

By Theorem \ref{thm:gaps}, the gap distribution of $(X,\omega)$ can be computed if we understand the return time function $R(x)$ of $h_s$ to $\Omega$ on each $\Omega_i$ in its local coordinates. 

We now describe some key properties of the return time function $R(x)$, many of which are explained in more detail in \cite{Kumanduri-Sanchez-Wang}. The first result we need is the following finiteness result. 

\begin{proposition}[\cite{Kumanduri-Sanchez-Wang}]
\label{prop:finiteness}
Each piece $\Omega_i$ of the transversal $\Omega$ breaks up into finitely many convex polygonal pieces. On each polygonal piece, there is a single vector $w_i = \begin{bmatrix} x \\ y \end{bmatrix}$ with $y > 0$ of $C_i(X,\omega)$ for which $M_{a,b}w_i$ has the smallest positive slope among all holonomy vectors with a horizontal length component of $\leq 1$. Then, the return time function at the point $(a,b) \in\Omega$ is given by $$R(a,b) = \frac{y}{a(ax+by)},$$ the slope of $M_{a,b}w_i$. On $\Omega_i$, $R(a,b)$ takes values in $(0,\infty)$. 
\end{proposition}

This idea of the return time function  is that $R(a,b)$ is the amount of time needed for the horocycle flow to pull the vector $M_{a,b}w_i$ down to become a short horizontal saddle connection so that $h_sM_{a,b}C_i (X,\omega)$ is in $\Omega$.


\subsection{Return time function bounds}
\label{sec:return-time-function-bounds}
We now wish to understand bounds on the return time function $R(a,b)$ on the transversal $\Omega$. The following two propositions will be necessary for the proofs of our main theorems. 

\begin{proposition} 
\label{prop:horocycle-return-time-bounds}
On each piece $\Omega_i$ of the transversal $\Omega$, the return time function $R(a,b)$ is bounded below by some constant $c_i > 0$. 
\end{proposition}
\begin{proof}
By Proposition \ref{prop:finiteness}, $\Omega_i$ breaks into finitely many polygonal pieces on which $R(a,b) = \frac{y}{a(ax+by)}$ for some $x \in \mathbb{R}$ and $y > 0$. On each of these pieces, the return time function is continuous and takes values in $(0,\infty)$. But $\overline{\Omega_i}$ is compact and so the denominator $a(ax+by)$ is bounded above by some $0 < M < \infty$ on the polygonal piece. Since $y > 0$ is fixed, $R(a,b)$ is then bounded below by $\frac{y}{M}> 0$. Taking $c_i$ to be the minimum of these bounds over all of the polygonal pieces, $R(a,b) > c_i > 0$ on all of $\Omega_i$
\end{proof}

\begin{proposition}
\label{prop:compactness-return-time-preimage}
Let $0 \leq c < d < \infty$. Then, $\{x \in \Omega : c \leq R(x) \leq d) \subset \Omega$ is compactly supported. 
\end{proposition}
\begin{proof}
It suffices to prove that $R^{-1}([c,d]) = \{x \in \Omega : c \leq R(x) \leq d) \subset \Omega$ is bounded away from the cusps of $\sltr/\Gamma$. As noted in Section \ref{sec:gaps-algorithm}, $\Omega$ breaks up into finitely many triangular pieces $\Omega_i$. In coordinates, a neighborhood of the vertex $(0,1/y_0)$ corresponds to the intersection of $\Omega_i$ with a neighborhood of a cusp of $\sltr/\Gamma$. 

By Proposition \ref{prop:finiteness}, each $\Omega_i$ breaks up into finitely many pieces on which the return time function is $R(a,b) = \frac{y}{a(ax+by)}$ with for some $x \in \mathbb{R}$ and $y > 0$. As $(a,b) \rightarrow (0,1/y_0)$, $R(a,b) \rightarrow \infty$. Thus, on each piece of each $\Omega_i$, $R^{-1}([c,d])$ is bounded away from $(0,1/y_0)$ and therefore bounded away from the cusps of $\sltr/\Gamma$. Thus, on $\Omega$, $R^{-1}([c,d])$ is bounded away from the cusps of $\sltr/\Gamma$ is therefore compactly supported. 
\end{proof}

\subsection{Proofs of main effective gaps theorems}
\label{sec:effective-gap-proofs}

Now that we have collected the necessary background, we will prove the main effective gap distribution theorems (Theorems \ref{thm:torus-gaps} and \ref{thm:Veech-gaps}) of this paper, as stated in the introduction. The proofs of both of these theorems follow directly from the corresponding effective equidistribution results of the horocycle flow on the transversal $\Omega$. 

\begin{proof}[Proof of Theorems \ref{thm:torus-gaps} and \ref{thm:Veech-gaps}]
Let $H_L$ be a length $L=1$ horocycle beginning at $(X,\omega)$ and $H_{L(t)} = g_t^{-1}H_L$ be a geodesic push of $H_L$ for $t=2 \log(R)$. 

Let $R: \Omega \rightarrow \mathbb{R}_+$ be the return time function of the horocycle flow $h_s$ to $\Omega$, and let $F$ be the distribution function of $R$. That is, let $$\int_a^b F(x) \, dx = m\{x\in \Omega : a \leq R(x) \leq b\}.$$ As shown in the proof of Theorem \ref{thm:gaps}, $\frac{|\mathbb{G}_R \cap (a,b)|}{N(R)}$ is the proportion of horocycle return times of $H_{L(t)} \cap \Omega$ that are between $a$ and $b$, and $\frac{|\mathbb{G}_R \cap (a,b)|}{N(R)} \rightarrow \int_a^b F(x) \, dx$ as $R\rightarrow \infty$. 

For any $0 \leq a < b < \infty$, let $f_{a,b}$ be the characteristic function on $R^{-1}([a,b]) \subset \Omega$. As shown in Proposition \ref{prop:compactness-return-time-preimage}, $f_{a,b}$ is compactly supported. $f_{a,b} \in L^2(\Omega)$ because it is a bounded function. Then, $\frac{|\mathbb{G}_R \cap (a,b)|}{N(R)} = \rho_{L(t)}(f_{a,b})$ for $t = 2 \log(R)$ and $\int_a^b F(x) \, dx = m(f_{a,b})$. The theorems then follows by applying Theorem \ref{thm:eff-equidistribution-on-transversal} for the torus and Theorem \ref{thm:eff-equidistribution-on-transversal-Veech-surfaces} for lattice surfaces and noting that $L(t) = R^2$. 
\end{proof}


\section{Horocycle length and intersection point bounds}
\label{sec:comparisons}

One of the main steps in the proof of an effective gaps theorem is to prove a result about effective equidistribution of long horocycles on the transversal $\Omega$. The main tool here is a thickening argument. One of the steps of this thickening argument is a comparison between the horocycle hitting measure $\rho_{L(t)}$ on $\Omega$ and the horocycle length measure $\nu_{L(t)}$ on $\sltr/\Gamma$, which is the topic of this section. 

\subsection{Spaces and measures}
\label{sec:spaces-measures}

As described in Section \ref{sec:gaps-algorithm}, for a lattice surface $(X,\omega)$, the transversal $\Omega \subset \sltr/\Gamma$ breaks up into finitely many pieces that are triangles in local $(a,b)$ coordinates. Then, $\Omega$ is a transversal to the horocycle flow and $\sltr/\Gamma$ can be thought of as a suspension over $\Omega$. Locally, over any piece of $\Omega$ with $(a,b)$ coordinates, $\sltr/\Gamma$ can be parametrized in $(a,b,s)$ coordinates where the $s$ comes from moving in $h_s$ direction. The height above each $(a,b)$ in this suspension is given by $R(a,b)$, the return time of the $h_s$ flow starting at $(a,b)$ in $\Omega$. 

When working in local $(a,b,s)$ coordinates, we sometimes refer to $\Delta$ as the \textbf{parameter space} for $\Omega$ and $S$ as the \textbf{suspension space} for $\sltr/\Gamma$. As defined in Equation \ref{eq:gt-push}, let $H_{L(t)}$ denote a $g_{t}^{-1}$ push of length $L$ horocycle segment based at $(X,\omega)$. We recall that $H_{L(t)}$ has length $L(t) = e^t\cdot L$, as commented upon in Remark \ref{rem:L(t)}. 

Then, on $\sltr/\Gamma$ we will work with the following two measures, both normalized so that the measure of $\sltr/\Gamma$ is $1$: 

\begin{enumerate}
    \item The Haar measure $\mu$. In local $(a,b,s)$ coordinates, 
\begin{equation}
    \label{eq:mu}
    \mu = c_\mu \, da \, db \, ds
\end{equation}
where the constant $c_\mu$ is chosen so that $\mu(\sltr/\Gamma) = \int_\Omega \int_0^{R(a,b)} c_\mu \, ds \, da \, db = 1$. 

\item The horocycle measure $\nu_{L(t)}$, which is the measure supported on the horocycle segment $H_L(t)$, normalized so that $\nu_{L(t)}(H_{L(t)}) = 1$.
\end{enumerate}

On the transversal $\Omega \subset \sltr/\Gamma$, we have the following two measures, also normalized so that the measure of $\Omega$ is $1$: 

\begin{enumerate}
    \item The Lebesgue measure $m$. In local coordinates, 
    \begin{equation}
    \label{eq:m}
    m = c_m \, da \, db,
\end{equation}
where $c_m$ is chosen so that $m(\Omega) = \int_\Omega c_m \, da \, db = 1$. 

    \item The horocycle counting measure $$\rho_{L(t)} = \frac{1}{\lvert H_{L(t)} \cap \Omega \rvert}\sum_{x \in H_{L(t)} \cap \Omega} \delta_{x}.$$
\end{enumerate}


\begin{remark} For the square torus, one can check (see \cite{Athreya-Cheung-poincare-section}) that $c_m=2$ and $c_\mu = \frac{1}{\zeta(2)}$. 
\end{remark}

\begin{remark}\label{rmk:smoothing} 
There is a map from the interior of the suspension space $S$ to the a subset of $SL_2(\R)$ (see Athreya-Cheung for the restriction of this map to the transversal~\cite{Athreya-Cheung-poincare-section}). The map is smooth, and consequently, a function that has been smoothed in the suspension space will pullback to a smooth function in $SL_2(\R)$.
\end{remark}



In Proposition \ref{prop:horocycle-return-time-bounds}, we saw that there exists a $w > 0$ such that the return time $R(a,b) > w$ for all points $(a,b)$ in every piece of $\Omega$. For the torus, the explicit return time function $R(a,b) = \frac{1}{ab}$ from Theorem \ref{thm:athreya-cheung} and the explicit parametrization of $\Omega$ given in Equation \ref{eq:poincare-section} give that we can take, for example, $w = \frac12$. 

Given a a measurable function $f: \Omega \to \R$, we can define a function $\tilde{f}$ that is thickened in the direction of the unstable horocycle flow $h_s$ as follows:
\begin{equation}\label{eq:tf}
    \tf (a,b,t) := \begin{dcases*}\tfrac{c_m}{c_\mu \cdot w} f(a,b)\chi_{[0,w]}(t) & if $(a,b) \in \Omega$, \\ 0 & otherwise\end{dcases*}.
\end{equation}


Here, $\chi_{[0,w]}$ denotes the indicator of the interval $[0,w]$. The constant multiplier in the definition of $\tf(a,b,t)$ is chosen so that the following proposition holds. 

\begin{proposition}
\label{prop:mu-and-m}
Given a measurable $f$, a $w$ so that the return time to $\Omega$ is bounded below by $w$, $m = c_m \, da \, db$ on $\Omega$ and $\mu = c_\mu \, da \, db \, ds$ on $\sltr/\Gamma$ in local $(a,b,s)$ coordinates, and a thickening $\tilde{f}$ as defined in Equation \ref{eq:tf}, $\mu(\tf) = m(f)$. 
\end{proposition}
\begin{proof}
    In local $(a,b)$ coordinates on the pieces of $\Omega$, we have that 
    \begin{align*}
    \mu(\tf)  = \int_\Omega \int_0^{R(a,b)} \frac{c_m}{c_\mu \cdot w} f(a,b) \chi_{[0,w]}(s) \, c_\mu \,ds \, da \, db = \frac{c_\mu \cdot w}{c_m} \cdot \frac{c_m}{c_\mu \cdot w} \int_\Omega f(a,b) \, c_m \, da \, db = m(f),
    \end{align*}
    where the factor of $w$ can be taken out of the integral because $R(a,b) > w$ for all $(a,b)$ in $\Omega$. 
\end{proof}

The goal of this section is to relate $\rho_{L(t)}(f)$ and $\nu_{L(t)}(\tf)$. To do so, we notice that 

$$\nu_{L(t)}(\tf) = \int_{\sltr/\Gamma} \frac{c_m}{c_\mu \cdot w} f(a,b) \chi_{[0,w]} d \nu_{L(t)} = \frac{c_m}{c_\mu \cdot w} \cdot w  \cdot \frac{1}{L(t)} \sum_{x \in H_{L(t)} \cap \Omega} f(x).$$
Since $\rho_{L(t)}(f) = \frac{1}{|H_{L(t)} \cap \Omega|} \sum_{x \in H_{L(t)} \cap \Omega} f(x)$, we have that 
\begin{equation}
    \label{eq:nu-and-rho}
    \nu_{L(t)}(\tilde{f}) = \frac{c_m}{c_\mu} \frac{|H_{L(t)} \cap \Omega|}{L(t)} \rho_{L(t)}(f).
\end{equation}

To understand the relationship between $\rho_{L(t)}(f)$ and $\nu_{L(t)}(\tf)$, we must therefore examine the relationship between $|H_{L(t)}\cap \Omega|$ and $L(t)$. We first look at the case of the square torus in Section \ref{sec:comparisons-torus} before moving on to general lattice surfaces in Section \ref{sec:comparisons-veech}. 

\subsection{Comparisons for the Square Torus}
\label{sec:comparisons-torus}


We begin with the case where $(X,\omega)$ is the area one square torus. In this case, the slopes of the saddle connections are rational and the slope gap distribution of the square torus is exactly the gap distribution of the Farey fractions. The gap distribution of the square torus was studied deeply by Athreya and Cheung in \cite{Athreya-Cheung-poincare-section}. 

In their paper, they consider the Farey fractions of denominator $\leq Q$, given by $$\mathcal{F}(Q) = \{\text{reduced fractions } \frac{p}{q} \text{ with } 0 \leq \frac{p}{q} \leq 1, q \leq Q\}.$$
When $(X,\omega)$ is the unit area square torus, these are exactly the slopes in $$\Lambda_Q = \{v = (a,b) \in \Lambda : 0 \leq b \leq a \leq Q\},$$ where $\Lambda$ are the saddle connection vectors of $(X,\omega)$. 

For the rest of our discussion on the torus, we will use the general notation of the pushed horocycles $H_{L(t)}$ rather than the $Q$ of Athreya-Cheung, which was specific to the square torus and Farey fractions. In Section \ref{sec:comparisons-veech}, we will use similar notation to do similar comparisons for general lattice surfaces. 

The Veech group $\Gamma$ of the square torus is $\Gamma = \sltz$. In this space, there is a family of closed horocycles corresponding to the single cusp of $X_2 = \sltr/\sltz$. We let $H_L$ be the closed horocycle of length $L=1$ in $\sltr/\sltz$. The following Theorem of Athreya and Cheung allows us to relate $|H_{L(t)} \cap \Omega|$ and $L(t) = e^t$, the length of $H_{L(t)}$. To relate this back to the notation of the paper of Athreya and Cheung, we have that $t = 2 \log(Q)$, $|H_{L(t)} \cap \Omega| = P(a,b)$, and $L(t)=s(a,b)$ where $(a,b) \in \Omega$ is $h_s$-periodic with period $Q^2$. 

\begin{remark}
The next few results are asymptotics in terms of $L(t)$. We note that since $L(t) = e^t$, these asymptotics could be written in terms of $t$ instead. However, we keep the $L(t)$ to make explicit the dependence on length. 
\end{remark}

\begin{theorem}[\cite{Athreya-Cheung-poincare-section} Theorem 1.5 therein]\label{thm:intersections-to-length}
Let $H_{L(t)}$ be the periodic horocycle in $\sltr/\sltz$ of period $L(t) = e^t$ . Then,
\begin{equation}
    |H_{L(t)} \cap \Omega| = N\left(\lfloor\sqrt{L(t)}\rfloor\right),
\end{equation}
 where $N(k)$ is the number of Farey fractions of denominator $\leq k$. 
\end{theorem}

Using this theorem, we have the following corollary, which is an intermediate step to relating $\rho_{L(t)}(f)$ with $\nu_{L(t)}(\tilde{f})$. 

\begin{corollary}\label{lem:torus-int-bound} Let $H_{L(t)}$ be the periodic horocycle in $\sltr/\sltz$ of length $L(t) = e^t$. Then, 
    \begin{equation} 
|H_{L(t)} \cap \Omega| = \frac{1}{2 \zeta(2)} L(t) + O\left(\sqrt{L(t)}\log(L(t))\right).
\end{equation}
\end{corollary}

\begin{proof}
    Let $\varphi(n)$ be Euler's totient function and $\Phi(k)$ be the summatory totient function. Since $\mathcal{F}(Q) = \mathcal{F}(Q-1) + \varphi(Q)$ and $\mathcal{F}(1) = |\{\frac01,\frac11\}| = 2$, it follows that 
    \begin{equation}
        N(k) = 1 + \sum_{n=1}^k \varphi(n) = 1 + \Phi(k).
    \end{equation} The summatory totient function is known (for example, see \cite{Hardy-Wright}) to have the asymptotic expansion 
     \begin{equation}
     \label{eq:Phi}
         \Phi(k) = \frac{1}{2\zeta(2)}k^2 + O(k \log k).
     \end{equation}

Putting this together, we have that 
\begin{align*} 
|H_{L(t)} \cap \Omega| & = N\left(\lfloor\sqrt{L(t)}\rfloor\right) = 1 + \Phi\left(\sqrt{L(t) + O(1)}\right) \\
& = \frac{1}{2 \zeta(2)} L(t) + O\left(\sqrt{L(t)}\log(L(t))\right).
\end{align*}
\end{proof}





\begin{lemma}\label{lem:rho-to-nu-Q}

Let $H_{L(t)}$ be the periodic horocycle on $\sltr/\sltz$ with length $L(t) = e^t$, let $f$ be a bounded, measurable function, and let $\tf$ be $f$ thickened by $w = \frac12$ in the horocycle flow direction. Then, there exists a constant, $K = K(f)$, independent of $t$, such that
\begin{equation*}
    |\rho_{L(t)} (f) - \nu_{L(t)} (\tf)| \leq K \lpar \frac{\log(L(t))}{\sqrt{L(t)}}\rpar
\end{equation*}

\end{lemma}


\begin{proof}
    We set $C:= \frac{1}{2 \zeta (2)}$ and $\err(t) := \sqrt{L(t)} \log(L(t))$. 

    Then, from Corollary \ref{lem:torus-int-bound}, there exists a constant such that, for $t$ sufficiently large,
\begin{align}\label{eq:T2-intersect-vs-length}
    |H_{L(t)} \cap \Omega| =  C\cdot L(t) +  O( \err(t)).
\end{align}

From Equation \ref{eq:nu-and-rho}, we have that $$\nu_{L(t)}(\tf) = \frac{c_m}{c_\mu} \frac{|H_{L(t)} \cap \Omega|}{L(t)} \rho_{L(t)}(f).$$ As remarked upon earlier, for the square torus, one can check (see \cite{Athreya-Cheung-poincare-section}) that $c_m=2$ and $c_\mu = \frac{1}{\zeta(2)}$. 

It follows then that $$\nu_{L(t)}(\tilde{f}) = \frac{c_m}{c_\mu} \cdot  \frac{C \cdot L(t) +  O( \err(t))}{L(t)} \rho_{L(t)}(f)= \rho_{L(t)}(f) \left(1 + O \left(\frac{c_m \cdot C}{c_\mu} \cdot \frac{\err(t)}{L(t)}\right)\right)  $$
\end{proof}

Let $M_f = \sup |f|$. Since $\rho_{L(t)}$ is normalized so that $\rho_L(t)(\Omega) = 1$, we have that $|\rho_{L(t)}(f)| \leq M_f$. Then, by substituting in $\err(t) = \sqrt{L(t)}\log(L(t))$, the statement of the Lemma follows for $K$ depending on $c_m$, $c_\mu$, $C$, $M_f$, and the implicit constant in the big-O term coming from Equation \ref{eq:Phi}. Thus, the constant $K$ depends on $f$ and is independent of $t$. 

\subsection{Comparisons for lattice surfaces}
\label{sec:comparisons-veech}




Following the convention of Athreya-Chaika-Leli\'{e}vre in \cite{Athreya-Chaika-Lelievre} we let $\Lambda_R(X,\omega)$ denote the set of saddle connections of slope between $0$ and $1$ for which the $x$ coordinate is positive and $\leq R$. 

As defined in Section \ref{sec:UW-proof}, let $H_L$ denote a horocycle of length $L$ beginning at $(X,\omega)$, thought of as a point in its $\sltr$ orbit $\sltr/\Gamma$. Let $H_{L(t)} = g_t^{-1}H_L$ be a push of $H_L$ by the geodesic flow, as defined in Equation \ref{eq:gt-push}.  Here, we will let $L=1$, so $L(t) = e^t$.


We wish to compare $L(t) = e^t$, with $|H_{L(t)} \cap \Omega|$, the number of times that $H_{L(t)}$ intersects the transversal $\Omega$. We note that this setup is the same as that of the square torus and Farey fraction discussion in Section \ref{sec:comparisons-torus}.  To estimate the latter quantity $|H_{L(t)} \cap \Omega|$, we will apply the following result of Burrin, Nevo, R\"{u}hr, and Weiss:

\begin{theorem}[\cite{Burrin}, Theorem 2.7 therein]
\label{thm:saddle-counts} 
Let $\Gamma < \sltr$ be a non-cocompact lattice, and let $Y :=\Gamma v$ be the $\Gamma$ orbit of a nonzero vector $v$. Let $S$ be a star-shaped domain with a non-negative piecewise Lipschitz boundary curve $\rho(\theta)$, and let $R \cdot S$ be the dilation of $S$ by $R > 0$. Then, for all $\epsilon > 0$,
\begin{equation}
    |Y \cap R\cdot S| = C_{Y,S} R^2 + O(R^{q_\Gamma + \veps}),    
\end{equation}
as $R \to \infty$, where the implied constant in the big-$O$ term depends on $\Gamma, Y,S,$ and $\veps$, and $q_\Gamma$ depends on the spectral gap of the automorphic representation of $\sltr$ on $L_0^2(\sltr/\Gamma)$. When $\Gamma$ is tempered, we can set $q_\Gamma = \frac74$. Otherwise, $q_{\Gamma} = 2 - \frac{1}{4n}$ for $n$ being the smallest even integer larger than $\frac{1}{1-s}$, where $\frac{1+s^2}{4}$ is the bottom of the spectrum of the hyperbolic Laplacian on $\Hh/\Gamma$ (non-inclusive of the 0). 
\end{theorem}

\begin{remark} 
\begin{enumerate}
\item This theorem is proven using a very general lattice counting argument in~\cite{Gorodnik-Nevo-erogidic-theory-of-lattice-subgroups} that depends on a mean ergodic theorem for the action of $\sltr$ on $\sltr/\gamma$. If one uses, for example, the Kunze-Stein phenomenon, or asymptotics of the Harish-Chandra function (see ~\cite{Gorodnik-Nevo-erogidic-theory-of-lattice-subgroups}, 5.2.2), one can confirm the values of $q_{\Gamma}$.
\item As this theorem will be instrumental to our rate bounds for effective gap distributions, any improvements made to the rate of convergence in Theorem \ref{thm:saddle-counts} would improve the effective gaps error rate of Theorem \ref{thm:Veech-gaps}. 
\end{enumerate}
\end{remark}

To relate $H_{L(t)}$ back to Theorem \ref{thm:saddle-counts}, we will let $t=2\log(R)$, so $L(t) = e^{2 \log(R)} = R^2$. We will adopt this convention for the rest of this section. As we will see in the proof of the following proposition, our choice of $R$ carries geometric significance beyond being the length of the horocycle.

\begin{proposition} Let $(X,\omega)$ be a lattice surface, $H_L$ a length $L=1$ horocycle beginning at $(X,\omega)$ in its orbit $\sltr/\Gamma$, and $H_{L(t)} = g_t^{-1}H_L$. Let $\epsilon > 0$ and $\mu = c_\mu \, ds \, da \, db$ and $m = c_m \, da \, db$ be measures on $\sltr/\Gamma$ and $\Omega$ respectively, both normalized to have total measure $1$. Then,  
\begin{equation}
    |H_{L(t)}\cap \Omega| = \frac{c_\mu}{c_m} (L(t)) + O(L(t)^{\frac12q_\Gamma + \veps}),
\end{equation}
as $t \to \infty$, where the implied constant in the big-$O$ term depends on $(X,\omega)$ and $\veps$, and $q_\Gamma$ is as in Theorem \ref{thm:saddle-counts}.
\end{proposition}
\begin{proof}

Let $t = 2 \log(R)$. By the definition of $\Omega$, horocycle $H_{L(t)}$ intersects the transversal $\Omega$ exactly when $h_s(g_{2\log(R)}^{-1}(X,\omega))$ has a short horizontal saddle connection. Since the action of $h_s$ on a given translation surface results in a vertical sheer and $L(t) = e^t = R^2$, the number of intersections is equal to the number of slopes of saddle connection vectors of $g_{2\log(R)}^{-1}(X,\omega)$ that reside in the triangle with vertices at $(0,0), (1,0),$ and $(1,R^2)$. Acting by $g_{2\log(R)}$, it is clear that this is equal to the number of slopes of saddle connections vectors of $(X,\omega)$ in the triangle with vertices at $(0,0), (R,0),$ and $(R,R)$. Thus, the number of intersections of the horocycle segment $H_{L(t)}$ with $\Omega$ is exactly the number of unique slopes of vectors in 
\begin{equation}
    \Lambda_R(X,\omega) = \{\text{saddle connections of } (X,\omega) \text{ with slope in } [0,1] \text{ and $x$ coordinate in } (0,1]\}.
\end{equation}
Figure \ref{fig:renormalization} depicts the effect of renormalization on the relevant saddle connections.

\begin{figure}[ht]
\centering
\begin{tikzpicture} 

\draw (-1,0) -- (3,0);
\draw (0,-1) -- (0,5);
\draw[fill=blue!25] (0,0) -- (2,0) -- (2,2)  -- (0,0);
\draw (4,0) -- (8,0);
\draw (5,-1) -- (5,5); 
\draw[fill = blue!25] (5,0) -- (6,0) -- (6,4) -- (5,0);
\draw[thick, ->] (3,2) -- (4,2);
\draw node at (3.5, 2.5) {$g_{-2 \log(R)}$};
\draw node at (2,-0.3) {$R$};
\draw node at (-0.3,2) {$R$};
\draw node at (6,-0.3) {$1$};
\draw node at (4.7, 4) {$R^2$};

\draw[dashed] (0,0) -- (1.8,1.6); 
\draw[dashed] (0,0) -- (1.6,0.6); 
\draw[dashed] (0,0) -- (1.4,0.2); 
\draw[dashed] (0,0) -- (0.8,0.6); 

\draw[fill] (1.8,1.6) circle [radius = 0.025]; 
\draw[fill] (1.6,0.6) circle [radius = 0.025]; 
\draw[fill] (1.4,0.2) circle [radius = 0.025]; 
\draw[fill] (0.8,0.6) circle [radius = 0.025]; 

\draw[dashed] (5,0) -- (5.9,3.2); 
\draw[dashed] (5,0) -- (5.8,1.2); 
\draw[dashed] (5,0) -- (5.7,0.4); 
\draw[dashed] (5,0) -- (5.4,1.2); 

\draw[fill] (5.9,3.2) circle [radius = 0.025]; 
\draw[fill] (5.8,1.2) circle [radius = 0.025]; 
\draw[fill] (5.7,0.4) circle [radius = 0.025]; 
\draw[fill] (5.4,1.2) circle [radius = 0.025]; 

\node at (1, -1.5) {$(X,\omega)$};
\node at (6, -1.5) {$g_{2\log(R)}^{-1}(X,\omega)$};
\end{tikzpicture}
\caption{The saddle connections $\Lambda_R$ on the surface $(X,\omega)$ on the left and the transformed saddle connections on $g_{2\log(R)}^{-1}(X,\omega)$ on the right.}
\label{fig:renormalization}
\end{figure}
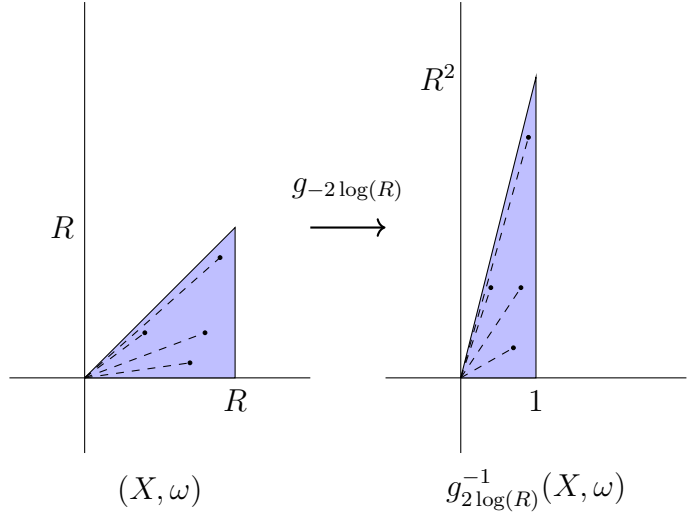

In Section \ref{sec:gaps-algorithm}, we discussed that the set of all saddle connection vectors of $(X,\omega)$ can be decomposed as 
\begin{equation}
\label{eq:gamma-orbits}
    \Lambda(X,\omega) = \bigcup_{i=1}^m (\Gamma v_i),
\end{equation}
the disjoint union of a finite number of $\Gamma$ orbits $\Gamma v_i$. 

Let $S$ be the triangle in $\mathbb{R}^2$ with vertices at $(0,0), (1,0),$ and $(1,1)$. For $t = 2 \log(R)$ wish to count $|H_{L(t)} \cap \Omega|$, which is equal to the number of unique slopes of $\Lambda_R (X,\omega) = \Lambda(X,\omega) \cap R\cdot S$. In Equation \ref{eq:gamma-orbits}, two saddle connections are in the same direction if and only if they are $\gamma v_i$ and $\gamma v_j$ for $v_i \neq v_j$ in the same direction (corresponding to the same cusp). Upon relabeling, we can take $v_1, \ldots, v_n$ to be a subset of the $v_i$'s that contains one vector for each cusp. Thus, we have that 
\begin{equation} 
\label{eq:gammaR-orbits}
|H_{L(t)} \cap \Omega| = |\text{unique slopes in } \Lambda_R(X,\omega)| = \left| \bigcup_{i=1}^n (\Gamma v_i) \cap (R\cdot S)\right|.
\end{equation}

We notice that $\bigcup_{i=1}^n (\Gamma v_i) \cap (R\cdot S) \subset \Lambda_R(X,\omega)$ contains exactly one saddle connection with each possible saddle connection slope direction in $\Lambda_R(X,\omega)$. 

 Then, Theorem \ref{thm:saddle-counts} gives us that  for each $v_i$,
 \begin{equation*}
     |\Gamma v_i \cap R\cdot S| = C_{\Gamma, v_i, S} R^2 + O(R^{q_\Gamma+\veps}).
 \end{equation*}
  By applying the decomposition from Equation \ref{eq:gammaR-orbits} and summing over the contributions of each $v_i$, we have that  
 \begin{equation} 
 \label{eq:int-counts}
 |H_{L(t)} \cap \Omega| = |\Lambda_R| = \sum_{i=1}^n |\Gamma v_i \cap R\cdot S| = C \cdot R^2 + O(R^{q_\Gamma + \veps}) = C\cdot L(t) + O(L(t)^{\frac12(q_\Gamma + \veps}) 
 \end{equation} 
 where the last equation follows since $L(t) = R^2$. Since $\veps$ was arbitrary, we can replace the $\frac12\veps$ with $\veps$ in the exponent of $L(t)$ in ther error term. We note that the $c$ depends on $\Gamma, S, v_1, \ldots, v_n$ and the constant in the $O$ term depends on $\Gamma, v_1, \ldots, v_n,$ and $\veps$. 

To determine $C$, we let $f$ be a compactly supported, smooth, $L^2$-integrable function with $m(f) > 0$. Then, by Proposition \ref{prop:mu-and-m}, $m(f) = \mu(\tf)$. By Equation \ref{eq:nu-and-rho}, $$\nu_{L(t)}(\tf) = \frac{c_m}{c_\mu} \cdot \frac{|H_{L(t)} \cap \Omega|}{L(t)} \rho_{L(t)}(f).$$ By the equidistribution of long horocycle segments (for example, see Theorem \ref{thm:eff-equid-long-hor-mixing}) and the proof of Theorem \ref{thm:gaps}, $\nu_{L(t)}(\tf) \rightarrow \mu(\tf)$ after approximating $\tf$ by a sequence of smooth functions, and $\rho_{L(t)}(f) \rightarrow m(f)$ as $t \rightarrow \infty$. For all of these statements to be consistent with Equation \ref{eq:int-counts}, it follows that $C = \frac{c_\mu}{c_m}$.
\end{proof}

We can then use this asymptotic to compare $|H_{L(t)} \cap \Omega|$ with $L(t) = R^2$.






\begin{lemma}\label{lem:rho-to-nu-Q-lattice}
Let $H_L$ be a horocycle of length $L$ in $\sltr/\Gamma$ and let $H_{L(t)} = g_t^{-1}H_L$ be a geodesic translate. Let $f$ be a bounded, measurable function. Then, there exists a constant, $K = K(f)$, independent of $t$, such that
\begin{equation*}
    |\rho_{L(t)} (f) - \nu_{L(t)} (\tilde{f})| \leq K L(t)^{-1 + \frac12 q_\Gamma + \varepsilon},
\end{equation*}

where $q_\Gamma$ is as in Theorem \ref{thm:saddle-counts}.

\end{lemma}
\begin{proof}
Set $C := c_\mu/c_m$ and $\err(t) := L(t)^{\frac12q_\Gamma + \veps}$. Using Theorem \ref{thm:saddle-counts} in place of Corollary \ref{lem:torus-int-bound} in the proof of Lemma \ref{lem:rho-to-nu-Q} gives the result by an identical argument.
\end{proof}

\begin{remark}
We note that since the square torus is a lattice surface where $\Gamma$ is tempered, we can apply Lemma \ref{lem:rho-to-nu-Q-lattice} and use $q_\Gamma = \frac74$ to get that $|\rho_{L(t)} (f) - \nu_{L(t)} (\tilde{f})| \leq K L(t)^{-\frac18 + \varepsilon}$. We notice that this bound is worse than the bound of $|\rho_{L(t)} (f) - \nu_{L(t)} (\tilde{f})| \leq K \log(L(t)) L(t)^{-\frac12}$ from Lemma \ref{lem:rho-to-nu-Q}, which was specific to the torus. 
\end{remark}

\section{Effective Estimates on Transversal}
\label{sec:effective-transversal}


In this section, we prove the main equidistribution theorems, Theorems \ref{thm:eff-equidistribution-on-transversal} and \ref{thm:eff-equidistribution-on-transversal-Veech-surfaces}. The proofs contain new notation, so we collect all of the notation in one place that be used in this section before stating the theorems. We begin with some matrix subgroups. 

\begin{align}
    A &:= \lbr g_t := \matr{e^{t/2}}{0}{0}{e^{-t/2}} : t \in \R \rbr\\
    N &:= \lbr h_s := \matr{1}{0}{-s}{1} : s \in \R \rbr \\
    U &:= \lbr u_s := \matr{1}{s}{0}{1} : s \in \R \rbr\\
    P &:= \lbr p_{a,b} := \matr{a}{b}{0}{a\inv} : a \in \R^*, b \in \R \rbr.
\end{align}

Here, the matrices in $A$ give the geodesic flow, $N$ the unstable horocycle flow, and $U$ the stable horocycle flow. The matrices in $P$ will be useful in defining the transversal $\Omega$ to the horocycle flow in local coordinates. Unless stated otherwise, in this paper, the ``horocycle flow" will refer to the unstable horocycle flow given by $h_s$. 

We recall that there are several relevant spaces and measures (see Section \ref{sec:spaces-measures} for the definitions). On $\sltr/\Gamma$, there is the Haar measure $\mu$ and the long horocycle measure $\nu_{L(t)}$ supported on the horocycle $H_{L(t)}$. On the transversal $\Omega \subset \sltr/\Gamma$, there is the Lebesgue measure $m$ and the horocycle counting measure $\rho_{L(t)}$ supported on $H_{L(t)} \cap \Omega$. All of these measures are normalized so that the measure of the whole space is $1$. 

We also define $S_r$ to be the parameter space associated with the suspension of $\Delta$: $$S_r := \{ (a,b,s) : , -r \leq s < R(a,b) - r \}.$$ Observe that Proposition \ref{prop:horocycle-return-time-bounds} implies that this definition makes sense for any $r$ less than the lower bound on the return time function. Moreover, on the torus, for any $r$, $\int_{S_r} 2dadbds = \int_{\Delta} R(a,b) 2dadb = 2 \zeta(2)$, where $2dadb$ is the Haar measure. See~\cite{Athreya-Cheung-poincare-section}. As in Section \ref{sec:comparisons}, we will use the notation $c_{\mu}$ for the constant such that for any $r$, $\int_{S_r} c_{\mu} dadbds = 1$.  Hence, for the torus, $c_{\mu} = \frac{1}{\zeta(2)}$. 


In what follows, we first prove the following equidistribution statement for the torus, first stated in the introduction. 



\effectivetorus*

\begin{remark} Theorem \ref{thm:eff-equidistribution-on-transversal} applies to periodic horocycles. We will use Lemma \ref{lem:rho-to-nu-Q} .
\end{remark}

One of the key inputs in the proof of the theorem is the effective equidistribution of long horocycle segments. We provide a description of this theorem in subsection \ref{sec:closed}. In subsection \ref{ssec:proof-equidistribution-on-section}, we set-up and prove Theorem \ref{thm:eff-equidistribution-on-transversal}. Each lemma we employ is completely general, with the exception of Lemma \ref{lem:rho-to-nu-Q}. In subsection \ref{sec:veech-transversal-equidistribution}, we replace Lemma \ref{lem:rho-to-nu-Q} with Lemma \ref{lem:rho-to-nu-Q-lattice} to prove the statement for general lattice surface, also first stated in the introduction. 

\effectiveveech*

\subsection{Effective equidistribution of long horocycles}
\label{sec:closed}
Let $\pi: SL_2(\R) \to L_0^2(SL_2(\R)/\Gamma)$ be the Koopman representation, where 
\begin{equation}\label{eq:Koopman-representation}
\pi(g)f(x) = f(g^{-1}x).
\end{equation}
Let $\omega = \begin{bmatrix} 0 & -1 \\ 1 & 0 \end{bmatrix}$ be a generator of the Lie algebra of $K$. A function $f \in L^2(SL_2(\R)/\Gamma, \mu)$ is called \emph{$K$-differentiable} if 
\begin{equation}\label{eq:K-derivative}
d(\omega) f := \lim_{h \to 0} \frac{ \pi(\mathrm{exp}(h \omega))f - f}{h}
\end{equation}
\noindent exists, where $h$ is a scalar and the convergence is in $L^2$. We say $f$ is a \emph{$K$-Sobolev function}, denoted $f \in \mathcal{S}_K$, if $f \in L^2(SL_2(\R)/\Gamma)$ and $d(\omega)f \in L^2(SL_2(\R)/\Gamma)$. We endow $\mathcal{S}_K$ with the norm 
\begin{equation}\label{eq:Sobolev-K-norm}
\mathcal{S}_K(f) = \left(\| f \|_2^2 + \| d(\omega)f \|_2^2 \right)^{\frac{1}{2}}\text{.}
\end{equation}

\begin{theorem}[Effective Equidistribution of Long Horocycles via Mixing]\label{thm:eff-equid-long-hor-mixing} 
Let $\Gamma$ be a non-cocompact lattice subgroup of $SL_2(\R)$ and $H_L$ be a the set of points corresponding to a segment of length $L$ of an (unstable) horocycle trajectory. Let $H_{L(t)} = g_t^{-1}H_L$, where $g_t$ is the geodesic flow for time $t$, and $L(t)$ denotes the length of $H_{L(t)}$. Let $f \in C_c^{\infty}(SL_2(\R)/\Gamma)$. 

\begin{equation*}
\left| \nu_{L(t)}(f) - \mu(f) \right| \leq \begin{cases} 
      C \mathcal{S}_{K}(f) \log(L(t)) L(t)^{-\frac{1}{10}} & \text{if } \Gamma \text{ is tempered}\\
      C(s) \mathcal{S}_{K}(f) L(t)^{-\frac{1}{10}(1-s)} & \text{if } \Gamma \text{ is non-tempered}
\end{cases} 
\end{equation*}

\noindent where $\frac{1-s^2}{4}$ is the bottom of the spectrum of the hyperbolic Laplacian on $\Hh/\Gamma$. 
\end{theorem}

\begin{remark}\label{rmk:long-horocycle-improvements} Sarnak was the first to give a full asymptotic expansion of the decay in~\cite{Sarnak}. We give a short proof of this well-known result in Appendix \ref{appendix:proof-of-effective-equidist} making the dependence on $f$ explicit. The proof provided uses the `Banana Trick', or `Margulis Thickening', pioneered by Margulis in his thesis~\cite{Margulis-thesis},~\cite{Wieser-note}. This technique, coupled with the effective mixing of the geodesic flow, yields the effective theorem. However, these exponents are not optimal, nor the best in the literature. Sarnak's result replaces the decay rate of $-\frac{1}{10}(1-s)$ with $-\frac{1}{2}(1-s)$ in the case of closed horocycles~\cite{Sarnak}. Similarly, for non-closed horocycles, we can replace the decay rate $-\frac{1}{10}(1-s)$ with $-\frac{1}{2}(1-s)$, provided we are willing to have test functions with at least four derivatives, and the corresponding Sobolev norm involving four derivatives. See Strombergsson~\cite{Strombergsson}.
\end{remark}

\subsection{Proof of equidistribution along the transversal}\label{ssec:proof-equidistribution-on-section}

To prove Theorem \ref{thm:eff-equidistribution-on-transversal}, we extend a compactly supported function $f \in L^2(\Omega)$ into the suspension space, smooth the extension, and employ effective equidistribution of long horocycles.  


Identify $L^2(\Omega, m)$ with $L^2(\Delta, c_m dadb)$ and work in the parameter space: let $f \in L^2(\Delta)$ be a compactly supported function. Trim the function away from the edges of $\Delta$: let $\mathcal{R_{\delta}}$ be the compact set in the parameter space where we trim $\delta < \frac{1}{2}R(a,b)$ (in Euclidean distance) off of each edge. In general, there may be finitely many components of the section $\mathcal{R}_{\delta}$, corresponding to the finitely many components of section $\Omega$. For the torus, we have 


\begin{equation*}
\psp_{\delta} = \{ (a,b) \in \R^2 : a,b \in (0,1-\delta], a+b \geq 1 + \delta \}.
\end{equation*}

Now define $$f_{\delta} = f \chi_{\psp_{\delta}} \text{.}$$ 

For the case of the torus, $c_m = 2$, and we can compute


\begin{equation}\label{eq:triangle-boundary-area}
\int_{\Delta} (f - f_{\delta})  \: \de a \de b \leq M_{\delta}\left(3\delta - \left(2 + \frac{1}{\tan(\pi/8)}\right)\delta^2\right)
\end{equation}

\noindent where $M_{\delta} = \sup_{x \in \Delta \setminus \psp_{\delta}}(f(x)) \leq M_f = \sup_{x \in \Delta} f(x)$.

For $0 < \varepsilon < \tfrac{\delta}{2} < \frac{1}{4} \inf R(a,b)$, define the \emph{open} set which is $\mathcal{R}_{\delta}$, but with an additional $\varepsilon$ trimmed off. For the torus, we have 
%
%
%
\begin{equation*}
\psp_{\delta + \veps} = \{ (a,b) \in \R^2 : a,b \in (0,1-\delta - \veps), a+b > 1 + \delta + \veps\} \text{.}
\end{equation*}

Next, we extend the function $f$ into the parameter space $S_{2\varepsilon}$ by extending the function in the unstable horocycle direction. We extend by the same method as in Section \ref{sec:comparisons}, with one exception: the suspension space that we extend the function into is $S_{2\varepsilon}$. However, this change is superficial. Recall that in $S_{2\varepsilon}$, $s \in [2\varepsilon, R(a,b) - 2\varepsilon)$. Define 

\begin{equation*}
    \tf (a,b,s) := \begin{dcases*} \frac{c_m}{c_{\mu} w} f(a,b)\chi_{[0,w]}(s) & if $(a,b) \in \Delta$, \\ 0 & otherwise\end{dcases*}.
\end{equation*}
 and  
\begin{equation*}
     \tf_{\delta}(a,b,s) = \begin{dcases*} \frac{c_m}{c_{\mu} w} f_\delta(a,b)\chi_{[0,w]}(s) & if $(a,b) \in \Delta$, \\ 0 & otherwise\end{dcases*}
\end{equation*}
where in both cases we must pick $w$ less than the lower bound on the return time function less $\varepsilon$. 

For the case of the torus, we choose $w = \frac{1}{2}$, and we have that $\frac{c_m}{c_{\mu} w} = 4\zeta(2)$ appears as the normalization factor. Note that since $\varepsilon < \frac{1}{2}$, and the minimum return time on the torus is $1$, $\tilde{f}$ and $\tilde{f}_{\delta}$ are well-defined with this choice of $w$.

In order to apply effective equidistribution of long horocycles to the trimmed function, we need to smooth this function. Following \ref{rmk:smoothing}, we opt to smooth the function in the suspension space. Define $\phi$ to be a bump function centered at $0$ in $\R^3$ such that 

\begin{equation*}
\phi(x) = \begin{cases} 
      k e^{\left({-\frac{1}{\left( 1 - \lvert x \rvert^2 \right)}}\right)} & \lvert x \rvert < 1 \\
      0 & \lvert x \rvert \geq 1
   \end{cases}
\end{equation*}

\noindent where $k$ is such that $\int_{\R^3} \phi \, dx = 1$. Then, for any $\varepsilon > 0$, define

\begin{equation*}
\phi_{\varepsilon}(x) = \frac{1}{\varepsilon^3} \phi\left(\frac{x}{\varepsilon}\right) .
\end{equation*}

\noindent We use this bump function to smooth $\tf_{\delta}$ via convolution. For $\varepsilon$ as above (the mollification is tied to the domains $\psp_{\delta}$ and $\psp_{\delta+\varepsilon}$ above), define
\begin{equation}
F_{\varepsilon, \delta}:= \phi_{\varepsilon}*\tilde{f}_{\delta} = \int \phi_{\varepsilon} (x - y) \tf_\delta  (x) \:\de x.  
\end{equation} 
where $x$ is a coordinate $(a,b,s) \in S_{2\varepsilon}$.

Fix $\varepsilon$ as above. Let $\rho_{L(t)}, \nu_{L(t)}, \mu$, and $m$ be as described at the beginning of Section \ref{sec:effective-transversal}. Let $f, \tilde{f}$, and $F_{\varepsilon, \delta}$ be as above. By the triangle inequality, we have the following.
\begin{align*}
\left| \rho_{L(t)}(f) - m(f) \right| &\leq \left| \rho_{L(t)}(f) - \nu_{L(t)}(\tilde{f}) \right| + \left| \nu_{L(t)}(\tilde{f}) - \nu_{L(t)}(F_{\varepsilon, \delta}) \right| + \left| \nu_{L(t)}(F_{\varepsilon, \delta}) - \mu(F_{\varepsilon, \delta}) \right| \\ 
&\,\,\,\,\, + \left| \mu(F_{\varepsilon, \delta}) - \mu(\tilde{f}) \right| + \left| \mu(\tilde{f}) - m(f) \right| 
\end{align*}

To prove Theorem \ref{thm:eff-equidistribution-on-transversal}, we must understand the rate at which each of these terms decays. Observe that the decay rate for the first term on the right hand side of the equation is given by Lemma \ref{lem:rho-to-nu-Q} in the case of the torus. Moreover, the last term is 0, as shown in Proposition \ref{prop:mu-and-m}. For the remainder of the section, we quantify the decay rates for each of the remaining terms.  

In the following lemma, we will need a parameter to track how far the support of a function extends into the cusp. We require the following definition: a \textbf{cusp region}, denoted $C(R)$ is a neighborhood of the cusp foliated by stable horocycles of measure $R$.

\begin{lemma}\label{lem:nu-smooth-to-nonsmooth} For sufficiently small $\varepsilon$, $\delta$, and $\tilde{\varepsilon}$, here exists a constant $B$ depending on $s$ and the geometry of the section such that $$\left\vert \nu_{L(t)}(\tilde{f}) - \nu_{L(t)}(F_{\varepsilon,\delta}) \right\vert \leq O\left( \varepsilon + \delta +\tilde{\varepsilon}^3\right) + B(s)\left(\frac{\varepsilon^{\frac{1}{2}}}{\tilde{\varepsilon}^2} + \frac{\delta^{\frac{1}{2}}}{\tilde{\varepsilon}^2}\right) L(t)^{-\frac{1}{10}(1-s)}.$$ Optimally, we have
$$\left\vert \nu_{L(t)}(\tilde{f}) - \nu_{L(t)}(F_{\varepsilon,\delta}) \right\vert \leq C(s) L(t)^{-\frac{3}{35}(1-s)}.$$ For the torus, $s=0$, and we have $$\left\vert \nu_{L(t)}(\tilde{f}) - \nu_{L(t)}(F_{\varepsilon,\delta}) \right\vert \leq C(s) L(t)^{-\frac{3}{35}}.$$
\end{lemma}

\begin{remark} The proof of Lemma \ref{lem:nu-smooth-to-nonsmooth} uses Theorem \ref{thm:eff-equid-long-hor-mixing}. If in lieu of the rates we have proven we used the results of Sarnak~\cite{Sarnak} for closed horocycles, then the $\frac{3}{35}$ can be replaced with $\frac{3}{7}$. For general horocycle segments, if we use the results of Strombergsson~\cite{Strombergsson}, we can replace the $\frac{3}{35}$ with $\frac{3}{19}$. The discrepancy between the closed and non-closed horocycles segments can be attributed to the fact that Str\"{o}mbergsson's results require Sobolev norms with four derivatives.  
\end{remark}

\begin{proof}
Observe that there exists a constant $C>0$ such that 
\begin{align*}
\left\vert \nu_{L(t)}(\tilde{f}) - \nu_{L(t)}(F_{\varepsilon,\delta}) \right\vert 
&\leq \frac{1}{L(t)} \int_0^{L(t)} \left\vert \tilde{f}(a(s),b(s),s) - F_{\varepsilon,\delta}(a(s),b(s),s) \right\vert ds \\ 
&\leq \frac{C}{L(t)} \int_0^{L(t)} \chi_{\overline{\Delta\setminus R_{\delta+\varepsilon}} \cap \mathrm{supp}(f) \times \left[-\varepsilon,w+\varepsilon\right]}(a(s),b(s),s) ds, 
\end{align*}
%
where $\chi_B$ indicates a characteristic function over a set $B$, since the functions in question only differ over small sets. Notice that $\overline{\Delta\setminus R_{\delta+\varepsilon}} \cap \mathrm{supp}(f) \times \left[-\varepsilon,w+\varepsilon\right]$ is a closed set, hence the characteristic function is defined over a closed set. We remark again that we can set $w= \frac{1}{2}$ for the torus.

To remove the dependence on the \emph{support} of $f$, we define an ancilliary parameter $C_f$, which is the supremum of the measures of cusp regions $C(R)$ such that $C(R) \cap \mathrm{supp(f)} = \emptyset$. Let $\hat{R}_{\delta + \varepsilon} = \overline{ R_{\delta + \varepsilon} \setminus \cup_j C_j(R) }$ where the union is taken over all cusp regions whose intersection with the support of $f$ is empty. In the case that the $\Gamma$ has multiple cusps, we run the argument for each cusp independently. Then, 
\begin{equation*}
\frac{C}{L(t)} \int_0^{L(t)} \chi_{\overline{\Delta\setminus R_{\delta+\varepsilon}} \cap \mathrm{supp}(f) \times \left[-\varepsilon,w+\varepsilon\right]}(a(s),b(s),s) ds \leq \frac{C}{L(t)} \int_0^{L(t)} \chi_{\overline{\Delta\setminus \hat{R}_{\delta+\varepsilon}} \times \left[-\varepsilon,w+\varepsilon\right]}(a(s),b(s),s) ds.
\end{equation*}
We extend this characteristic function and apply Theorem \ref{thm:eff-equid-long-hor-mixing}. To do this, we opt to work directly in $SL_2(\R)$ to avoid boundary considerations in the suspension space. 

Let $K$ be the image of $\overline{\Delta\setminus \hat{R}_{\delta+\varepsilon}} \times \left[-\varepsilon,w+\varepsilon\right]$ in $SL_2(\R)/\Gamma$. Define a neighborhood of the identity $\mathcal{O}_{\tilde{\varepsilon}} \subset SL_2(\R)$ for any $\tilde{\varepsilon} > 0$ as follows: 
\begin{equation}\label{eq:O-tilde-epsilon}
\mathcal{O}_{\tilde{\varepsilon}} := \{  u_{s_1}g_{s_2} h_{s_3} : s_1, s_3 \in (-\tilde{\varepsilon}, \tilde{\varepsilon}), s_2 \in \left[ 0 ,\ln(1+ \tilde{\varepsilon}/2) \right) \}.
\end{equation}
Fix $\tilde{\varepsilon}$ to be smaller than the injectivity radius of the compact set $K$. Let $K_{\tilde{\varepsilon}} = \mathcal{O}_{\tilde{\varepsilon}}K$. Observe that 
\begin{equation*}
\frac{C}{L(t)} \int_0^{L(t)} \chi_{\overline{\Delta\setminus \hat{R}_{\delta+\varepsilon}} \times \left[-\varepsilon,w+\varepsilon\right]}(a(s),b(s),s) ds \leq \frac{C}{L(t)} \int_0^{L(t)} \chi_{K_{\tilde{\varepsilon}}}(g\Gamma(s)) ds. 
\end{equation*}
 
Recall that $\left\{ H, X^+, X^- \right\}$ forms a basis for $\mathfrak{sl}_2(\R)$, the Lie algebra of $SL_2(\R)$ where 

\begin{equation*}
H = \begin{bmatrix} 1 & 0 \\ 0 & -1 \end{bmatrix}, X^+ = \begin{bmatrix} 0 & 1 \\ 0 & 0 \end{bmatrix}, X^- = \begin{bmatrix} 0 & 0 \\ 1 & 0 \end{bmatrix}.
\end{equation*}

With this basis, define an inner product $\langle x, y \rangle = x^t y$ for $x,y \in \mathfrak{sl}_2(\R)$. Recall that $\mathfrak{sl}_2(\R)$ is identified with the tangent space at the identity in $SL_2(\R)$, so we have the induced inner product on this tangent space. Using left multiplication on the group, we can pullback the metric to any tangent space to generate a smooth left-invariant metric $g$ on $SL_2(\R)$. Take the volume form of this metric, and observe that it is bi-invariant, hence a multiple of the Haar measure. Furthermore, for $x \in T_p(SL_2(\R))$, define $\vert x \vert_p = \sqrt{g_p( x, x )}$, and let $d(p_1, p_2)$ denote the induced distance between $p_1$ and $p_2$ for $p_1, p_2 \in SL_2(\R)$.   

Define $\phi: SL_2(\R) \to \R$:

\begin{equation*}
\phi(g) = \begin{cases} 
      k_1 e^{\left({-\frac{1}{\left( k_2^2 - d(g, e)^2 \right)}}\right)} & d(g, \mathrm{Id}) < k_2 \\
      0 & d(g, \mathrm{Id}) \geq k_2
   \end{cases}
\end{equation*}

\noindent where $k_1$ is such that $\int_{G} \phi = 1$ and $k_2$ is such that there exists a diffeomorphism between $B_{k_2}(\mathrm{Id})$ and a small neighborhood of $0$ in the Lie algebra $\mathfrak{sl}_2(\R)$.
Now define 
\begin{equation*}
\phi_{\tilde{\varepsilon}}(g) = \frac{1}{\tilde{\varepsilon}^3} \phi\left({\mathrm{exp}\left(\frac{\log(g)}{\tilde{\varepsilon} \lvert \log(g) \rvert_{\mathrm{Id}}}\right)}\right)
\end{equation*}
and observe that since $\tilde{\varepsilon}$ is less than the injectivity radius of $K$, 
\begin{align*}
\frac{C}{L(t)} \int_0^{L(t)} \chi_{\overline{\Delta\setminus \hat{R}_{\delta+\varepsilon}} \times \left[-\varepsilon,w+\varepsilon\right]}(a(s),b(s),s) ds &\leq \frac{C}{L(t)} \int_0^{L(t)} \phi_{\tilde{\varepsilon}}*\chi_{K_{\tilde{\varepsilon}}}(g\Gamma(s)) ds \\
&\leq \frac{C}{L(t)} \int_0^{L(t)} \chi_{K_{\tilde{\varepsilon}}}(g\Gamma(s)) ds.
\end{align*}
Apply Theorem \ref{thm:eff-equid-long-hor-mixing} to the middle integral.  
\begin{align*}
\left\vert \frac{C}{L(t)} \int_0^{L(t)} \phi_{\tilde{\varepsilon}}*\chi_{K_{\tilde{\varepsilon}}}(g\Gamma(s)) ds - C \int_{\sltr/\Gamma} \phi_{\tilde{\varepsilon}}*\chi_{K_{\tilde{\varepsilon}}}(g\Gamma) d\mu(g\Gamma) \right\vert &\leq  C(s) \mathcal{S}_{K}(\phi_{\tilde{\varepsilon}}*\chi_{K_{\tilde{\varepsilon}}}) L(t)^{-\frac{1}{10}(1-s)}.
\end{align*}
Hence, we have
\begin{align*}
\left\vert \nu_{L(t)}(\tilde{f}) - \nu_{L(t)}(F_{\varepsilon,\delta}) \right\vert &\leq \frac{C}{L(t)} \int_0^{L(t)} \phi_{\tilde{\varepsilon}}*\chi_{K_{\tilde{\varepsilon}}}(g\Gamma(s)) ds \\
&\leq C \int_{\sltr/\Gamma} \phi_{\tilde{\varepsilon}}*\chi_{K_{\tilde{\varepsilon}}}(g\Gamma) d\mu(g\Gamma) +  C(s) \mathcal{S}_{K}(\phi_{\tilde{\varepsilon}}*\chi_{K_{\tilde{\varepsilon}}}) L(t)^{-\frac{1}{10}(1-s)}. \\
\end{align*}
Recall that $\mathrm{supp}(\phi_{\tilde{\varepsilon}} * \chi_{K_{\tilde{\varepsilon}}}) \subset \overline{\mathrm{supp}(\phi_{\varepsilon}) + \mathrm{supp}(\chi_{K_{\tilde{\varepsilon}}})}$. The measure of the support of $\phi_{\tilde{\varepsilon}}$ is proportional to $\tilde{\varepsilon}^3$. We claim that there exists a $B$ such that $\mu(K_{\tilde{\varepsilon}}) \leq (1 + B\tilde{\varepsilon}) \mu(K)$, and refrain from proving this until after completing the rest of the proof. 

The measure of the support of $K_{\tilde{\varepsilon}} = \mathcal{O}_{\tilde{\varepsilon}} K$ is proportional to $\delta + \varepsilon$ and higher degree terms. In the case of the torus, an explicit computation shows that the measure is proportional to $$(1+B\tilde{\varepsilon})\left(w + 2 \varepsilon \right)\left(3\left(\delta+ \varepsilon \right) - \left(2+ \frac{1}{\tan(\frac{\pi}{8})}\right)(\delta + \varepsilon)^2\right),$$ where the measure of $K$ is estimated using Equation \ref{eq:triangle-boundary-area}. Consequently, the first term in the previous inequality is proportional to $\tilde{\varepsilon}^3 
 + \delta + \varepsilon$ along with other higher degree terms.

We estimate the Sobolev norm as follows (in the case of the torus).
\begin{align*}
\mathcal{S}_K(\phi_{\tilde{\varepsilon}}*\chi_{K_{\tilde{\varepsilon}}})^2 &= \| \phi_{\tilde{\varepsilon}}* \chi_{K_{\tilde{\varepsilon}}} \|_2^2 + \| d(\omega) \phi_{\tilde{\varepsilon}}*\chi_{K_{\tilde{\varepsilon}}} \|_2^2 \\
& \leq \| \chi_{K_{\tilde{\varepsilon}}} \|_2^2 + \frac{1}{\tilde{\varepsilon}^4} \cdot \sup_g \left(d\omega( \phi) \right) \| \chi_{K_{\tilde{\varepsilon}}} \|_2^2 \\
& \leq (1+B\tilde{\varepsilon})\left(w + 2 \varepsilon \right)\left(3\left(\delta+ \varepsilon \right) - \left(2+ \frac{1}{\tan(\frac{\pi}{8})}\right)(\delta + \varepsilon)^2\right)  \left(1 + \frac{1}{\tilde{\varepsilon}^4} \cdot \sup_g \left(d\omega( \phi) \right)\right). \\
\end{align*}
%
%
For the general case, the orders of $\tilde{\varepsilon},\varepsilon$, and $\delta$ remains the same.  

Combining terms, for any $0 < \varepsilon < \frac{\delta}{2} < \frac{1}{4} \inf R(a,b)$, and for any sufficiently small $\tilde{\varepsilon} > 0$, we set $\varepsilon = L(t)^{-\alpha}$, $\delta = L(t)^{-\beta}$, and $\tilde{\varepsilon} = L(t)^{-\eta}$ for any $\alpha$, $\beta$, $\eta > 0$, and deduce that the optimal rate is achieved by choosing $\alpha = \beta = 3\eta = \frac{3}{35}(1-s)$. Then, we have that here exists a constant $C(s)$ such that 
\begin{align*}
 \left\vert \nu_{L(t)}(\tilde{f}) - \nu_{L(t)}(F_{\varepsilon,\delta}) \right\vert \leq C(s) L(t)^{-\frac{3}{35}(1-s)}.
\end{align*}
In the case of the torus, since $SL_2(\Z)$ is tempered, we can set $s=0$.
\end{proof}

To complete the proof of Lemma \ref{lem:nu-smooth-to-nonsmooth} we must prove the claim. We remark that the claim is straightforward: our set $K$ is the image of a set defined in $UAN$-coordinates in the suspension space. Our perturbation of this set is exactly a small perturbation in each of the coordinate directions. For sufficiently small perturbations (sufficiently small $\tilde{\varepsilon}$), the main increase in the volume is of order $\tilde{\varepsilon}$. However, there is one technicality we must address, which is that small perturbations in the $(a,b,s)$-coordinate will push a $(1,b,s)$-coordinate outside of the suspension space. There are two remedies we could choose from: first, we could understand how the edges of the suspension space are identified, or second, we could pull the set back into $SL_2(\R)$, change coordinates, and do the computation there. We will use the latter approach.

\begin{lemma}\label{lem:well-roundedness} Let $K$ be a compact set of $SL_2(\R)/\Gamma$, and let $\mathcal{O}_{\tilde{\varepsilon}}$ be as in Equation \ref{eq:O-tilde-epsilon}. Let $K_{\tilde{\varepsilon}} \leq \mathcal{O}_{\tilde{\varepsilon}} K$. For $\tilde{\varepsilon}$ sufficiently small, there exists a constant $C>0$ independent of $\tilde{\varepsilon}$ such that $$\mu(K_{\tilde{\varepsilon}}) = (1 + C\tilde{\varepsilon}) m(K).$$
\end{lemma}

\begin{proof} $K$ is the image of the set $\overline{\Delta\setminus \hat{R}_{\delta+\varepsilon}} \times \left[-\varepsilon,\frac{1}{2}+\varepsilon\right]$ in $SL_2(\R)/\Gamma$ in $SL_2(\R)$. Observe that points of the form $(1,b,s)$ may be in the set, which lies on the boundary of the suspension space. To avoid complications coming from working with elements on the boundary of $S$, we will continue working in $SL_2(\R)/\Gamma$. Recall that for any $g \in \Delta\setminus \hat{R}_{\delta+\varepsilon} \times \left[-\varepsilon,\frac{1}{2}+\varepsilon\right]$ in $SL_2(\R)/\Gamma$, we can write this in the form $p_{a,b}h_s$. By applying elements from $\mathcal{O}_{\varepsilon}$, we see how the parameters $a,b,$ and $s$ change. For sufficiently small $\tilde{\varepsilon}$, there are 6 extreme points in the set $\overline{\mathcal{O}_{\tilde{\varepsilon}}}$, each corresponding to a choice of $\pm \tilde{\varepsilon}$ for each of the three parameters $s_1, s_2,$ and $s_3$. Applying each of these elements to a point $p_{a,b}h_s \in K$ gives us a maximum distortion of the parameters $a,b$ and $s$. For instance, applying the element $u_{\tilde{\varepsilon}}g_{\tilde{\varepsilon}}h_{\tilde{\varepsilon}} \in \overline{\mathcal{O}_{\tilde{\varepsilon}}}$, we see 
\begin{align*}
a &\longrightarrow \left(\frac{2 + \tilde{\varepsilon}}{2-2\tilde{\varepsilon}ab} \right) \cdot a  \\
b &\longrightarrow \left( \frac{2+ \tilde{\varepsilon}}{2} - \tilde{\varepsilon}^2\left( \frac{2}{2+\tilde{\varepsilon}} \right) \right)\cdot b + \tilde{\varepsilon} \left( \frac{2}{2+\tilde{\varepsilon}} \right)\cdot a \\
s &\longrightarrow \left(\frac{2}{2+ \tilde{\varepsilon}} \right) \left(\frac{2 + \tilde{\varepsilon}}{2-2\tilde{\varepsilon}ab} \right) \cdot s + \tilde{\varepsilon}\left( \frac{2}{2+ \tilde{\varepsilon}} \right) \left( \frac{2 + \tilde{\varepsilon}}{2-2\tilde{\varepsilon}ab} \right) \cdot a(a-sb).
\end{align*}
To understand the additional measure we pick up when perturbing the set $K$, we re-write coordinates in terms of the $UAK$-decomposition (this is called the $NAK$-decomposition in other literature, but we have reserved $N$ for the unstable horocycle). For $g = p_{a,b}h_s$, we write $g = u(g)\tilde{a}(g)k(g)$ for $u(g) = u_u= \begin{bmatrix} 1 & u \\ 0 & 1  \end{bmatrix}$, (and apologize for the redundancy in the use of the $u$ variable), $\tilde{a}(g) = g_{t}$, and $k(g) = k_{\theta} = \begin{bmatrix} \cos(\theta) & \sin(\theta) \\ -\sin(\theta) & \cos(\theta)  \end{bmatrix}$ and compute
\begin{align*}
\theta(a,b,s) &= \mathrm{arccot}(s) \\
t(a,b,s) &= 2\log \left(\frac{1}{\sqrt{1 + s^{-2}}}\right) \\
u(a,b,s) &= ab - \frac{a^2s^{-1}}{1 + s^{-2}}.
\end{align*}
There appear to be singularities in these functions when $s = 0$, but we can compute: $t \to -\infty$ as $s \to 0^+$, $t \to \infty$ as $s \to 0^-$, and $u = ab$ when $s=0$. The angle $\theta$ has a discontinuity here. For $s<0$,  $\theta$ will approach $-\pi/2$ as $s$ approaches $0$. On the other hand, for $s>0$, $\theta$ approaches $\pi/2$ as $s$ approaches $0$. However, this does not hinder our computation: we will only be concerned with the extreme values of $s$ in the set $K$, and here, $s \neq 0$.

By applying  $u_{\tilde{\varepsilon}}g_{\tilde{\varepsilon}}h_{\tilde{\varepsilon}} \in \overline{\mathcal{O}_{\tilde{\varepsilon}}}$ to an element in $K$, and looking at the $UAK$-decomposition, we obtain perturbed coordinates $(\tilde{u}, \tilde{t}, \tilde{\theta})$, and by appling $u_{\tilde{\varepsilon}_1}$, $g_{\tilde{\varepsilon}_2}$, and $r_{\tilde{\varepsilon}_3}$ to an element in $K$: 
\begin{align}
\tilde{\theta}(u, t, \theta) &= \arctan(-\frac{c}{d}) \\
\tilde{t}(u, t, \theta) &= \log(d^2-c^2) \\
\tilde{s}(u, t, \theta) &= \frac{a + d(c^2-d^2)}{c}
\end{align}
where 
\begin{align*}
a &= \left( \frac{2}{2+\tilde{\varepsilon}_2}- \tilde{\varepsilon}_1\tilde{\varepsilon}_3 \frac{2}{2+\tilde{\varepsilon}}\right)\left( \cos(\theta)e^{\frac{t}{2}} - u \sin(\theta)e^{-\frac{t}{2}} \right) + \left(\tilde{\varepsilon}_1\frac{2}{2+\tilde{\varepsilon}_2} \right) \left( -\sin(\theta)e^{-\frac{t}{2}} \right)\\
b &= \left( \frac{2}{2+\tilde{\varepsilon}_2}- \tilde{\varepsilon}_1\tilde{\varepsilon}_3 \frac{2}{2+\tilde{\varepsilon}} \right)\left( \sin(\theta) e^{\frac{t}{2}} + u \cos(\theta)e^{-\frac{t}{2}}\right) +  \left( \tilde{\varepsilon}_1\frac{2}{2+\tilde{\varepsilon}_2} \right) \left( \cos(\theta) e^{-\frac{t}{2}} \right)\\
c &= \left(-\tilde{\varepsilon}_3\frac{2}{2+\tilde{\varepsilon}_2} \right)\left( \cos(\theta)e^{\frac{t}{2}} - u \sin(\theta)e^{-\frac{t}{2}} \right) +\left( \frac{2}{2+\tilde{\varepsilon}_2}\right)\left( -\sin(\theta)e^{-\frac{t}{2}} \right)\\
d &= \left( -\tilde{\varepsilon}_3\frac{2}{2+\tilde{\varepsilon}_2}  \right)\left( \sin(\theta) e^{\frac{t}{2}} + u \cos(\theta)e^{-\frac{t}{2}} \right)+\left(\frac{2}{2+\tilde{\varepsilon}_2} \right)\left( \cos(\theta)e^{-\frac{t}{2}} \right)\\
\end{align*}
Observe that $ad-bc = 1$. 

Thus, the perturbation is of order $\tilde{\varepsilon}_i$ in each coordinate. By letting $\tilde{\varepsilon}_i = \pm \tilde{\varepsilon}$ for each $i \in \{1, 2, 3 \}$, and integrating both $K$ and $K_{\tilde{\varepsilon}}$ over $SL_2(\R)$, we observe that 
$\frac{\mu(K_{\tilde{\varepsilon}})}{\mu(K)} = 1 + C_1\tilde{\varepsilon} + C_2\tilde{\varepsilon}^2 + C_3 \tilde{\varepsilon}^3$, for some constants $C_1, C_2,$ and $C_3$, as desired. 
\end{proof}

\begin{lemma}\label{lem:mu-smooth-to-nonsmooth} 
For $\varepsilon, \delta > 0$, we have that  
$$\left\vert \mu(F_{\varepsilon, \delta}) - \mu(\tilde{f}) \right\vert = O(\delta + \varepsilon \delta)$$ where the implied constant depends on $f$. In the case of the torus, we have $$\left\vert \mu(F_{\varepsilon, \delta}) - \mu(\tilde{f}) \right\vert \leq C (w + 2\varepsilon) \left(3\delta - \left(2 + \frac{1}{\tan(\pi/8)}\right)\delta^2\right).$$
\end{lemma}

\begin{proof}
We estimate $\tilde{f}$ using the function $\tilde{f}_{\delta}$, and apply the triangle inequality.  

\begin{align*}
\labs \mu(F_{\eta,\delta}) - \mu(\tilde{f}) \rabs &\leq \labs \mu(F_{\eta, \delta}) - \mu(\tilde{f}_{\delta}) \rabs + \labs \mu(\tilde{f}_{\delta}) - \mu(\tilde{f}) \rabs \\
&\leq \left\vert \int_{S_{2\varepsilon}}  \phi_{\varepsilon} * \tilde{f}_{\delta} \, \de \mu - \int_{S_{2\varepsilon}} \tilde{f}_{\delta} \, \de \mu \right\vert + \int_{S_{2\varepsilon}} \left\vert \tilde{f}_{\delta} - \tilde{f} \right\vert \de \mu\\
&= \left\vert  \int_{\R^3} \phi_{\varepsilon}\, \de \mu \int_{S_{2\varepsilon}} \tilde{f}_{\delta} \, \de \mu - \int_{S_{2\varepsilon}} \tilde{f}_{\delta} \, \de \mu \right\vert  + \int_{S_{2\varepsilon}} \left\vert \tilde{f}_{\delta} - \tilde{f} \right\vert \de \mu \\
&= \int_{S_{2\varepsilon}} \left\vert \tilde{f}_{\delta} - \tilde{f} \right\vert \de \mu 
\end{align*}

The functions $\tilde{f}_\delta$ and $\tilde{f}$ only differ on $\Delta \setminus R_{\delta} \cap \mathrm{supp}(f)$. On the torus, by Equation \ref{eq:triangle-boundary-area}, we have 
\begin{align*}
\int_{S_{2\varepsilon}} \left\vert \tilde{f}_{\delta} - \tilde{f} \right\vert d\mu \leq \frac{1}{\zeta(2)} M_{\tilde{f}}(w + 2\varepsilon) \left(3\delta - \left(2 + \frac{1}{\tan(\pi/8)}\right)\delta^2\right)
\end{align*}
where $M_{\tilde{f}}$ is the essential supremum of the function $\tilde{f}$. In the general case, the orders of $\varepsilon$ and $\delta$ remain the same. 
\end{proof}

The main theorem follows from these lemmas. 

\begin{proof}[Proof of Theorem \ref{thm:eff-equidistribution-on-transversal}]
Recall that 
\begin{align*}
\left| \rho_{L(t)}(f) - m(f) \right| &\leq \left| \rho_{L(t)}(f) - \nu_{L(t)}(\tilde{f}) \right| + \left| \nu_{L(t)}(\tilde{f}) - \nu_{L(t)}(F_{\varepsilon}) \right| + \left| \nu_{L(t)}(F_{\varepsilon}) - \mu(F_{\varepsilon}) \right| \\ 
&\,\,\,\,\, + \left| \mu(F_{\varepsilon}) - \mu(\tilde{f}) \right| + \left| \mu(\tilde{f}) - m(f) \right| 
\end{align*}
Apply Lemma \ref{lem:rho-to-nu-Q}, Lemma \ref{lem:nu-smooth-to-nonsmooth}, Lemma \ref{lem:mu-smooth-to-nonsmooth}, and Proposition \ref{prop:mu-and-m} to the first, second, fourth and fifth terms, respectively. Observe that Lemma \ref{lem:rho-to-nu-Q} only applies to \emph{closed} horocycle trajectories. For the third term, apply Theorem \ref{thm:eff-equid-long-hor-mixing}: for any $\varepsilon$ and $\delta$ sufficiently small,
\begin{equation*}
\left| \nu_{L(t)}(F_{\varepsilon,\delta}) - \mu(F_{\varepsilon,\delta}) \right| \leq \begin{cases} 
      C \mathcal{S}_{K}(F_{\varepsilon,\delta}) \log(L(t)) L(t)^{-\frac{1}{10}} & \text{if } \Gamma \text{ is tempered}\\
      C(s) \mathcal{S}_{K}(F_{\varepsilon,\delta}) L(t)^{-\frac{1}{10}(1-s)} & \text{if } \Gamma \text{ is non-tempered.}
\end{cases} 
\end{equation*}
We can estimate the Sobolev norm. Assume $\varepsilon$ is sufficiently small so that the measure of the support of $\phi_{\varepsilon}$ is less than the injectivity radius of the compact support of $\tilde{f}_{\delta}$. There exists a constants $B, D>0$ such that
\begin{align*}
\mathcal{S}_K(F_{\varepsilon,\delta})^2 &= \mathcal{S}_K(\phi_{\varepsilon}*\tilde{f}_{\delta})^2 \\
&= \| \phi_{\varepsilon} * \tilde{f}_{\delta} \|_2^2 + \| d(\omega) \phi_{\varepsilon} * \tilde{f}_{\delta} \|_2^2 \\
& \leq \| \tilde{f}_{\delta} \|_2^2 + \frac{1}{\varepsilon^4} \cdot \sup_g \left(d\omega( \phi) \right) \| \tilde{f}_{\delta} \|_2^2 \\
& \leq D \| f_{\delta} \|_2^2 \left(1+ \frac{B}{\varepsilon^4}\right) \\
& \leq D \| f \|_2^2 \left(1+ \frac{B}{\varepsilon^4}\right).
\end{align*}
Combining all of the terms, for any $0 < \varepsilon < \frac{\delta}{2} < \frac{1}{4} \inf R(a,b)$ where $\varepsilon$ is sufficiently small, and for any sufficiently small $\tilde{\varepsilon} > 0$, we set $\varepsilon = L(t)^{-\alpha}$, $\delta = L(t)^{-\beta}$, and $\varepsilon{\delta} = L(t)^{-\eta}$ for any $\alpha$, $\beta$, $\eta > 0$, and deduce that the optimal rate is achieved by choosing $\alpha = \beta = \frac{1}{30}(1-s)$. Note that this corresponds to choosing $\eta = \frac{1}{6}(1-s)$

\begin{equation*}
\left| \rho_{L(t)}(f) - m(f) \right| \leq \begin{cases} 
      C \| f \|_2 \log(L(t)) L(t)^{-\frac{1}{30}} & \text{if } \Gamma \text{ is tempered}\\
      C(s) \| f \|_2 L(t)^{-\frac{1}{30}(1-s)} & \text{if } \Gamma \text{ is non-tempered.}
\end{cases}
\end{equation*}
In the case of the torus, $SL_2(\Z)$ is tempered, so we set $s=0$.
\end{proof}

\begin{remark}\label{rmk:rates-for-torus} The proof of Theorem \ref{thm:eff-equidistribution-on-transversal} uses Theorem \ref{thm:eff-equid-long-hor-mixing} in multiple places. If in lieu of the rates we have proven we used the results of Sarnak~\cite{Sarnak} for closed horocycles, then the $\frac{1}{30}$ can be replaced with $\frac{1}{6}$. In the case of the torus, since we are using Lemma \ref{lem:rho-to-nu-Q}, the result is not valid for general horocycle segments.  
\end{remark}

\begin{remark} In the proof, we could avoid ``clipping" the function $f$ (to create $f_{\delta}$) by working in $SL_2(\R)/\Gamma$ instead of the parameter space $S_{2\varepsilon}$. This would eliminate the $\delta$ parameter, however, since the order of $\delta$ matches $\varepsilon$ throughout the proof, with the exception of the $\frac{1}{\varepsilon^2}$ appearing because of the Sobolev norm, we cannot improve the final decay rate by removing $\delta$. The parameters are $\varepsilon$ and $\tilde{\varepsilon}$ establish the rate independent of $\delta$.
\end{remark}

\subsection{General lattice surfaces: equidistribution along the transversal}
\label{sec:veech-transversal-equidistribution}

To prove Theorem \ref{thm:eff-equidistribution-on-transversal-Veech-surfaces}, we follow the same technique as in Section \ref{ssec:proof-equidistribution-on-section}: we extend a compactly supported function $f \in L^2(\Omega)$ into the suspension space, smooth the extension, and employ effective equidistribution of long horocycles. Note that in the case of a general Veech surfrace, there may be more than one cusp. Here, we fix one cusp, and observe that the rate is the same for all. 

As before, identify $L^2(\Omega)$ with $L^2(\Delta)$. Lemmas \ref{lem:nu-smooth-to-nonsmooth}, \ref{lem:mu-smooth-to-nonsmooth}, and \ref{prop:mu-and-m} apply to our setting. We remark that in Lemma \ref{lem:nu-smooth-to-nonsmooth}, we rely on Equation \ref{eq:triangle-boundary-area} to estimate the measure of the set $K$. In the general case, Equation \ref{eq:triangle-boundary-area} does not take that form, but the orders of $\varepsilon$ and $\delta$ remain the same. 

In lieu of Lemma \ref{lem:rho-to-nu-Q}, we use Lemma \ref{lem:rho-to-nu-Q-lattice}. We observe that the rate of the general case, as applied to the torus, is slower than the rate for the torus in Theorem \ref{thm:eff-equidistribution-on-transversal}. We attribute this to our choice to use the lattice point counting result \ref{thm:saddle-counts} to prove Lemma \ref{lem:rho-to-nu-Q-lattice}. While Theorem \ref{thm:saddle-counts} is both an effective and remarkable result, one may be able to improve it in our setting. 

\begin{proof}[Proof of Theorem \ref{thm:eff-equidistribution-on-transversal-Veech-surfaces}]
Recall that 
\begin{align*}
\left| \rho_{L(t)}(f) - m(f) \right| &\leq \left| \rho_{L(t)}(f) - \nu_{L(t)}(\tilde{f}) \right| + \left| \nu_{L(t)}(\tilde{f}) - \nu_{L(t)}(F_{\varepsilon}) \right| + \left| \nu_{L(t)}(F_{\varepsilon}) - \mu(F_{\varepsilon}) \right| \\ 
&\,\,\,\,\, + \left| \mu(F_{\varepsilon}) - \mu(\tilde{f}) \right| + \left| \mu(\tilde{f}) - m(f) \right| 
\end{align*}
Apply Lemma \ref{lem:rho-to-nu-Q-lattice}, Lemma \ref{lem:nu-smooth-to-nonsmooth}, Lemma \ref{lem:mu-smooth-to-nonsmooth}, and Proposition \ref{prop:mu-and-m} to the first, second, fourth and fifth terms, respectively. Observe that Lemma \ref{lem:rho-to-nu-Q-lattice} applies to both closed and non-closed horocycle trajectores. For the third term, apply Theorem \ref{thm:eff-equid-long-hor-mixing} exactly as above. Moreover, the estimate of the Sobolev norm is the same as above. 

The difference between Theorem \ref{thm:eff-equidistribution-on-transversal} and Theorem \ref{thm:eff-equidistribution-on-transversal-Veech-surfaces} lies in the balancing of the decay rates, but only if we use the best possible rates in the literature (see Remark \ref{rmk:rates-for-Veech-surfaces} below). With our set-up, with Lemma \ref{lem:rho-to-nu-Q-lattice}, the slower rate still comes from balancing the out terms coming from the Sobolev norms. We observe
\begin{equation*}
\left| \rho_{L(t)}(f) - m(f) \right| \leq \begin{cases} 
      C \| f \|_2 \log(L(t)) L(t)^{-\frac{1}{30}} & \text{if } \Gamma \text{ is tempered}\\
      C(s) \| f \|_2 L(t)^{-\frac{1}{30}(1-s)} & \text{if } \Gamma \text{ is non-tempered,}
\end{cases}
\end{equation*}
as desired.
\end{proof}

\begin{remark}\label{rmk:rates-for-Veech-surfaces} 
\begin{enumerate} 
\item The proof of Theorem \ref{thm:eff-equidistribution-on-transversal-Veech-surfaces} uses Theorem \ref{thm:eff-equid-long-hor-mixing} in multiple places, just as Theorem \ref{thm:eff-equidistribution-on-transversal} does. If in lieu of the rates we have proven we used the results of Sarnak~\cite{Sarnak} for closed horocycles, then the $\frac{1}{30}$ can be replaced with $\frac{1}{8n} - \varepsilon$ for any $\varepsilon>0$ where $n$ is the smallest even number larger than $\frac{1}{1-s}$. However, Sarnak's result makes the balanced terms \emph{faster} than the contribution coming from the rates in Lemma \ref{lem:rho-to-nu-Q-lattice}, and these rates comes from Theorem \ref{thm:saddle-counts}. To improve our result for lattice surfaces for closed horocycles, one would need to improve the lattice counting result first, and likely our methods, second. 
\item For general horocycle segments, if we use the results of Strombergsson~\cite{Strombergsson}, we can replace the $\frac{1}{30}$ with $\frac{1}{18}$, just as in the case of the torus. The discrepancy between the closed and non-closed horocycles segments can be attributed to the fact that Stromberggson's results require Sobolev norms with four derivatives. 
\end{enumerate}
\end{remark}

These effective equidistribution results for the torus and for general lattice surfaces (Theorems \ref{thm:eff-equidistribution-on-transversal} and Theorems \ref{thm:eff-equidistribution-on-transversal-Veech-surfaces}) are the key input needed for the proofs of the main effective gaps results of this paper (Theorem \ref{thm:torus-gaps} and Theorem \ref{thm:torus-gaps}). The effective gaps theorems were proven earlier in the paper in Section \ref{sec:effective-gap-proofs}.


\bibliographystyle{alpha}
\bibliography{gaps-bibliography}

\newcommand{\etalchar}[1]{$^{#1}$}
\begin{thebibliography}{BMMM{\etalchar{+}}23}

\bibitem[AC12]{Atheya-Chaika}
Jayadev~S. Athreya and Jon Chaika.
\newblock The distribution of gaps for saddle connection directions.
\newblock {\em Geom. Funct. Anal.}, 22(6):1491--1516, 2012.

\bibitem[AC14]{Athreya-Cheung-poincare-section}
Jayadev~S. Athreya and Yitwah Cheung.
\newblock A {P}oincar\'{e} section for the horocycle flow on the space of lattices.
\newblock {\em Int. Math. Res. Not. IMRN}, (10):2643--2690, 2014.

\bibitem[ACL15]{Athreya-Chaika-Lelievre}
Jayadev~S. Athreya, Jon Chaika, and Samuel Leli\`evre.
\newblock The gap distribution of slopes on the golden {L}.
\newblock In {\em Recent trends in ergodic theory and dynamical systems}, volume 631 of {\em Contemp. Math.}, pages 47--62. Amer. Math. Soc., Providence, RI, 2015.

\bibitem[ACZ15]{Athreya-Cobeli-Zaharescu}
Jayadev~S. Athreya, Cristian Cobeli, and Alexandru Zaharescu.
\newblock Radial density in {A}pollonian packings.
\newblock {\em Int. Math. Res. Not. IMRN}, (20):9991--10011, 2015.

\bibitem[BMMM{\etalchar{+}}23]{Berman-2n}
Jonah Berman, Taylor McAdam, Ananth Miller-Murthy, Caglar Uyanik, and Hamilton Wan.
\newblock Slope gap distribution of saddle connections on the {$2n$}-gon.
\newblock {\em Discrete Contin. Dyn. Syst.}, 43(1):1--56, 2023.

\bibitem[BNRW20]{Burrin}
Claire Burrin, Amos Nevo, Rene R\"{u}hr, and Barak Weiss.
\newblock Effective counting for discrete lattice orbits in the plane via {E}isenstein series.
\newblock {\em Enseign. Math.}, 66(3-4):259--304, 2020.

\bibitem[BV16]{Browning-Vinogradov}
Tim Browning and Ilya Vinogradov.
\newblock Effective {R}atner theorem for {$\text{SL}(2,\mathbb{R})\ltimes\mathbb{R}^2$} and gaps in {$\sqrt{n}$} modulo 1.
\newblock {\em J. Lond. Math. Soc. (2)}, 94(1):61--84, 2016.

\bibitem[Edw21]{Edwards}
Samuel~C. Edwards.
\newblock On the rate of equidistribution of expanding translates of horospheres in {$\Gamma\backslash G$}.
\newblock {\em Comment. Math. Helv.}, 96(2):275--337, 2021.

\bibitem[GN10]{Gorodnik-Nevo-erogidic-theory-of-lattice-subgroups}
Alexander Gorodnik and Amos Nevo.
\newblock {\em The Ergodic Theory of Lattice Subgroups}.
\newblock Princeton University Press, Princeton, 2010.

\bibitem[HS06]{Hubert-Schmidt}
Pascal Hubert and Thomas~A. Schmidt.
\newblock An introduction to {V}eech surfaces.
\newblock In {\em Handbook of dynamical systems. {V}ol. 1{B}}, pages 501--526. Elsevier B. V., Amsterdam, 2006.

\bibitem[HT92]{Howe-Tan}
Roger Howe and Eng~Chye Tan.
\newblock {\em Non-{Abelian} harmonic analysis. {Applications} of {{\(SL(2,{\mathbb{R}})\)}}}.
\newblock Universitext. New York etc.: Springer-Verlag, 1992.

\bibitem[HW08]{Hardy-Wright}
G.~H. Hardy and E.~M. Wright.
\newblock {\em An introduction to the theory of numbers}.
\newblock Oxford University Press, Oxford, sixth edition, 2008.
\newblock Revised by D. R. Heath-Brown and J. H. Silverman, With a foreword by Andrew Wiles.

\bibitem[KM12]{Kleinbock-Margulis-09}
D.~Y. Kleinbock and G.~A. Margulis.
\newblock On effective equidistribution of expanding translates of certain orbits in the space of lattices.
\newblock In {\em Number theory, analysis and geometry}, pages 385--396. Springer, New York, 2012.

\bibitem[KSW24]{Kumanduri-Sanchez-Wang}
Luis Kumanduri, Anthony Sanchez, and Jane Wang.
\newblock Slope gap distributions of veech surfaces.
\newblock {\em Algebraic and Geometric Topology}, 24(2):951--980, April 2024.

\bibitem[Li15]{Li-Frobenius}
Han Li.
\newblock Effective limit distribution of the {F}robenius numbers.
\newblock {\em Compos. Math.}, 151(5):898--916, 2015.

\bibitem[LMW22]{LMW-Effective}
Elon Lindenstrauss, Amir Mohammadi, and Zhiren Wang.
\newblock Effective equidistribution for some one parameter unipotent flows, 2022.

\bibitem[LMW23]{LMW-flat-torus}
Elon Lindenstrauss, Amir Mohammadi, and Zhiren Wang.
\newblock Quantitative equidistribution and the local statistics of the spectrum of a flat torus.
\newblock {\em J. Anal. Math.}, 151(1):181--234, 2023.

\bibitem[Mar04]{Margulis-thesis}
G.~A. Margulis.
\newblock {\em On some aspects of the theory of {Anosov} systems. {With} a survey by {Richard} {Sharp}: {Periodic} orbits of hyperbolic flows. {Transl}. from the {Russian} by {S}. {V}. {Vladimirovna}}.
\newblock Springer Monogr. Math. Berlin: Springer, 2004.

\bibitem[Mar10]{Marklof-frobenius}
Jens Marklof.
\newblock The asymptotic distribution of {F}robenius numbers.
\newblock {\em Invent. Math.}, 181(1):179--207, 2010.

\bibitem[Mas88]{Masur88}
Howard Masur.
\newblock Lower bounds for the number of saddle connections and closed trajectories of a quadratic differential.
\newblock In {\em Holomorphic functions and moduli, {V}ol. {I} ({B}erkeley, {CA}, 1986)}, volume~10 of {\em Math. Sci. Res. Inst. Publ.}, pages 215--228. Springer, New York, 1988.

\bibitem[Mas90]{Masur90}
Howard Masur.
\newblock The growth rate of trajectories of a quadratic differential.
\newblock {\em Ergodic Theory Dynam. Systems}, 10(1):151--176, 1990.

\bibitem[NRW20]{Nevo-Ruhr-Weiss}
Amos Nevo, Rene R\"{u}hr, and Barak Weiss.
\newblock Effective counting on translation surfaces.
\newblock {\em Advances in Mathematics}, 369, 2020.

\bibitem[San22]{Sanchez}
Anthony Sanchez.
\newblock Gaps of saddle connection directions for some branched covers of tori.
\newblock {\em Ergodic Theory Dynam. Systems}, 42(10):3191--3245, 2022.

\bibitem[Sar81]{Sarnak}
Peter Sarnak.
\newblock Asymptotic behavior of periodic orbits of the horocycle flow and {E}isenstein series.
\newblock {\em Comm. Pure Appl. Math.}, 34(6):719--739, 1981.

\bibitem[Str13]{Strombergsson}
Andreas Str\"{o}mbergsson.
\newblock On the deviation of ergodic averages for horocycle flows.
\newblock {\em J. Mod. Dyn.}, 7(2):291--328, 2013.

\bibitem[Str15]{strombergsson-effective-lift}
Andreas Str\"ombergsson.
\newblock An effective {R}atner equidistribution result for {$\mathrm{SL}(2,\mathbb{R})\ltimes\mathbb{R}^2$}.
\newblock {\em Duke Math. J.}, 164(5):843--902, 2015.

\bibitem[UW15]{Uyanik-Work}
Caglar Uyanik and Grace Work.
\newblock {The Distribution of Gaps for Saddle Connections on the Octagon}.
\newblock {\em International Mathematics Research Notices}, 2016(18):5569--5602, 11 2015.

\bibitem[Vor05]{Vorobets}
Yaroslav Vorobets.
\newblock Periodic geodesics on generic translation surfaces.
\newblock In {\em Algebraic and topological dynamics}, volume 385 of {\em Contemp. Math.}, pages 205--258. Amer. Math. Soc., Providence, RI, 2005.

\bibitem[Wie]{Wieser-note}
Andreas Wieser.
\newblock On the ``banana"-trick of margulis.

\bibitem[Wri15]{Wright}
Alex Wright.
\newblock Translation surfaces and their orbit closures: an introduction for a broad audience.
\newblock {\em EMS Surv. Math. Sci.}, 2(1):63--108, 2015.

\bibitem[Zor06]{Zorich}
Anton Zorich.
\newblock Flat surfaces.
\newblock In {\em Frontiers in number theory, physics, and geometry. {I}}, pages 437--583. Springer, Berlin, 2006.

\end{thebibliography}
\appendix
\section{Effective Equidistribution}\label{appendix:proof-of-effective-equidist}
In what follows, we provide a proof of Theorem \ref{thm:eff-equid-long-hor-mixing}. We begin by stating the effective mixing for the geodesic flow on quotients of $\sltr/\Gamma$. 

\begin{theorem}[\cite{Nevo-Ruhr-Weiss} Theorem 3.2 therein, \cite{Howe-Tan} Proposition 3.1.5 therein]\label{thm:decay-of-matrix-coeff} 
Let $SL_2(\R)$ act ergodically on $(SL_2(\R)/\Gamma, \mu)$. There exists $C>0$ such that for any $f_1, f_2 \in L^2_0(SL_2(\R)/\Gamma), \mu$ which is $K$-differentiable and $t>0$, we have  

\begin{equation*}
\left\vert \langle \pi(g_t) f_1 , f_2 \rangle \right\vert \leq
\begin{cases} 
      C t e^{-\frac{t}{2}} \mathcal{S}_{K}(f_1)\mathcal{S}_{K}(f_2) & \text{if } \Gamma \text{ is tempered}\\
      C(s) e^{-\frac{t}{2}(1-s)} \mathcal{S}_{K}(f_1)\mathcal{S}_{K}(f_2) & \text{if } \Gamma \text{ is non-tempered}
   \end{cases} 
\end{equation*}

\noindent where $\frac{1-s^2}{4}$ is the bottom of the spectrum of the Laplacian on $\Hh/\Gamma$.
\end{theorem}

\begin{remark}
Note that in the non-tempered case, the constant $C$ depends on $s$. Since we will fix a surface first, we do not require a $C$ independent of $s$. However, it is still possible to have a result with a $C$ independent of $s$, but the quality of the estimate declines. See the proof of Theorem 3.2 in~\cite{Nevo-Ruhr-Weiss}. For the values of the constants $C$ and $C(s)$, see Proposition 3.1.5 in~\cite{Howe-Tan}.
\end{remark}

Before providing the proof of Theorem \ref{thm:eff-equid-long-hor-mixing}, we make two observations. First, let $H_{L} =\{h_s h\Gamma : s \in [0,L] \}$ denote a segment of an unstable horocycle of length $L$. Observe that $g_t^{-1}\cdot H_{L} = H_{L(t)}(g_t^{-1} h\Gamma) := \{ h_{s} g_t^{-1}h\Gamma : 0 \leq s \leq L(t) \}$ where 

\begin{proposition} $L(t) = e^t L \text{.}$
\end{proposition}

\begin{proof} This follows from the fact that $g_t^{-1} h_{s} = h_{se^t}g_t^{-1}$.
\end{proof}

Second, we can smooth functions in $SL_2(\R)$. As in Lemma \ref{lem:nu-smooth-to-nonsmooth}, recall that $\left\{ H, X^+, X^- \right\}$ forms a basis for $\mathfrak{sl}_2(\R)$, the Lie algebra of $SL_2(\R)$ where 

\begin{equation*}
H = \begin{bmatrix} 1 & 0 \\ 0 & -1 \end{bmatrix}, X^+ = \begin{bmatrix} 0 & 1 \\ 0 & 0 \end{bmatrix}, X^- = \begin{bmatrix} 0 & 0 \\ 1 & 0 \end{bmatrix}.
\end{equation*}

With this basis, define an inner product $\langle x, y \rangle = x^t y$ for $x,y \in \mathfrak{sl}_2(\R)$. Recall that $\mathfrak{sl}_2(\R)$ is identified with the tangent space at the identity in $SL_2(\R)$, so we have the induced inner product on this tangent space. Using left multiplication on the group, we can pullback the metric to any tangent space to generate a smooth left-invariant metric $g$ on $SL_2(\R)$. Take the volume form of this metric, and observe that it is bi-invariant, hence a multiple of the Haar measure. Furthermore, for $x \in T_p(SL_2(\R))$, define $\vert x \vert_p = \sqrt{g_p( x, x )}$, and let $d(p_1, p_2)$ denote the induced distance between $p_1$ and $p_2$ for $p_1, p_2 \in SL_2(\R)$.   

Define $\phi: SL_2(\R) \to \R$:

\begin{equation*}
\phi(g) = \begin{cases} 
      k_1 e^{\left({-\frac{1}{\left( k_2^2 - d(g, e)^2 \right)}}\right)} & d(g, \mathrm{Id}) < k_2 \\
      0 & d(g, \mathrm{Id}) \geq k_2
   \end{cases}
\end{equation*}

\noindent where $k_1$ is such that $\int_{G} \phi = 1$ and $k_2$ is such that there exists a diffeomorphism between $B_{k_2}(\mathrm{Id})$ and a small neighborhood of $0$ in the Lie algebra $\mathfrak{sl}_2(\R)$.
For $\varepsilon > 0$, define 
\begin{equation*}
\phi_{\varepsilon}(g) = \frac{1}{\varepsilon^3} \phi\left({\mathrm{exp}\left(\frac{\log(g)}{\varepsilon \lvert \log(g) \rvert_{\mathrm{Id}}}\right)}\right).
\end{equation*}

\noindent We will convolve $\phi_{\varepsilon}$ with functions to smooth them in $SL_2(\R)$. 

Recall the definition of convolution. Let $G$ be a locally compact second countable group with unimodular Haar measure $\mu$. Let $f_1: G \to \R$ and $f_2: G \to \R$ be locally integrable functions. The \textbf{convolution}, denoted $f_1*f_2$, is 
\begin{align*}
f_1*f_2(g) :&= \int_G f_1(h)f_2(h^{-1}g) d\mu(h) = \int_G f_1(gh) f_2(h^{-1}) d\mu(h) \\
&= \int_G f_1(h)f_2(h^{-1}g) d\mu(h) = \int_G f_1(gh^{-1})f_2(h) d\mu(h),
\end{align*}
where the second line follows from the fact that $G$ is unimodular.
\begin{proof}[Proof of Theorem \ref{thm:eff-equid-long-hor-mixing}]   
Let $\mathcal{O}_P$ a neighborhood of the identity in the subgroup of upper triangular matrices $P = UA \subset SL_2(\R)$ such that $$g_t \mathcal{O}_P g_t^{-1} \subset \mathcal{O}_P.$$ It is well-known that such a neighborhood exists. See, for instance~\cite{Wieser-note}. 

Let $m_P^l$ denote the left Haar measure on $P$ and let $m_{N}$ denote the right Haar measure on $N$. Then any right Haar measure on $SL_2(\R)$ restricted to $PN$ is proportional to the pushforward $\phi_{*}(m_P^l \times m_N)$ where $\phi: P \times N \to SL_2(\R)$, $\phi(p_{a,b},n_s) = p_{a,b}n_s$.

Without loss of generality, assume that the unstable horocycle segment $S_{L}$ has no self-intersections, otherwise pull it back via $g_t^{-1}$. 

Let $\delta_0$ be the injectivity radius of $S_{L(0)}$ (a compact segment of an unstable horocycle). Fix $\mathcal{O}_P \subset P \cap B_{\delta_0}(\mathrm{Id})$ where $B_{\delta_0}(\mathrm{Id}) \subset SL_2(\R)$. Let $m_P^l(\mathcal{O}_P) = \delta$. We will shrink this set, taking $\delta \to 0$, at a rate specified at the end of the proof. 

Now let $f \in C_c^{\infty}(SL_2(\R)/\Gamma)$. Note that $f$ is uniformly continuous. We have that the integral of interest can be approximated by an integral in $SL_2(\R)$ by thickening the horocycle using $\mathcal{O}_P$. Indeed, by the uniform continuity of the $f$, for any $\tilde{\varepsilon} > 0$, we can pick $\mathcal{O}_P$ sufficiently small so that   
\begin{align*} 
&\left\vert \frac{1}{e^t L} \int_{0}^{e^t L} f(h_{s} g_t^{-1} h\Gamma) \, ds - \frac{1}{e^t L} \int_{0}^{e^t L}  \left( \frac{1}{m_P^l(g_t \mathcal{O}_P g_t^{-1})} \int_{g_t \mathcal{O}_P g_t^{-1}} f(b h_{s} g_t^{-1} h\Gamma) dm_P^l \right)\, ds \right\vert \\
&\leq   \frac{1}{e^t L} \int_{0}^{e^t L}  \left( \frac{1}{m_P^l(g_t \mathcal{O}_P g_t^{-1})} \int_{g_t \mathcal{O}_P g_t^{-1}}\left\vert f(h_{s}g_t^{-1} h\Gamma) - f(bh_{s}g_t^{-1} h\Gamma) \right\vert  \right)\ \, dm_P^l ds \\
&< \tilde{\varepsilon} \\
&= O(\sqrt{\delta})
\text{.}
\end{align*}

To approximate the integral over the thickened set, we observe that it can be written as an integral corresponding to the mixing properties of the geodesic flow on $SL_2(\R)/\Gamma$. First, let $S_t := \{ bh_{s} : 0 \leq s < e^t L, \, b \in g_t \mathcal{O}_P g_t^{-1}  \}$. Then 
\begin{align*}
\frac{1}{e^t L} \int_{0}^{e^t L}  \left( \frac{1}{m_P^l(g_t \mathcal{O}_P g_t^{-1})} \int_{g_t \mathcal{O}_P g_t^{-1}} f(b h_{s} g_t^{-1} h\Gamma) dm_P^l \right)\, ds & = \frac{1}{m_{SL_2(\R)}(S_t)} \int_{S_t} f(g g_t^{-1} h\Gamma) dm_{SL_2(\R)}(g) \\
\end{align*}
\noindent where we are using the aforementioned decomposition of the Haar measure. Now, substitute $g$ with $g_t^{-1} g g_t$. Since $SL_2(\R)$ is unimodular, the measure is unchanged. Further, observe that $g_t^{-1}S_tg_t = S_0$. 
\begin{equation*}
\frac{1}{m_{SL_2(\R)}(S_t)} \int_{S_t} f(g g_t^{-1} h\Gamma) dm_{SL_2(\R)}(g) = \frac{1}{m_{SL_2(\R)}(S_0)} \int_{S_0} f(g_t^{-1} g(h\Gamma)) \, dm_{SL_2(\R)}(g) \\
\end{equation*}
Now, observe that the projection of $S_0$ into $SL_2(\R)/\Gamma$ is injective, provided we pick $\mathcal{O}_P$ sufficiently small. Let $P_0 = \{g (h\Gamma) : g \in S_0 \}$ be the projection. Then,
\begin{equation*}
\frac{1}{m_{SL_2(\R)}(S_0)} \int_{S_0} f(g_t^{-1} g(h\Gamma)) \, dm_{SL_2(\R)}(g) = \frac{1}{\mu(P_0)} \int_{P_0} f(g_t^{-1} g(h\Gamma)) \, d\mu(g (h\Gamma)) \\
\end{equation*}
Now, observe that the integral is a matrix coefficient. As $t \to \infty$, we expect exponential decay. 
\begin{align*}
\frac{1}{\mu(P_0)} \int_{P_0} f(g_t^{-1} g(h\Gamma)) \, d\mu(g (h\Gamma)) &= \int_{SL_2(\R)/\Gamma} f(g_t^{-1} g(h\Gamma)) \frac{\chi_{P_0}(g(h\Gamma))}{\mu(P_0)} \, d\mu(g (h\Gamma)) \\
&= \left\langle \pi(g_t) f, \frac{\chi_{P_0}}{\mu(P_0)} \right\rangle 
\end{align*}
Let $\chi$ denote $\frac{\chi_{P_0}}{\mu(P_0)}$. We cannot directly apply Theorem \ref{thm:decay-of-matrix-coeff} since $\chi$ is not $K$-smooth. Let $\psi$ be a continuous approximation of $\chi$ where for any small $\varepsilon>0$, $\int \vert \chi - \psi \vert \leq \frac{\varepsilon}{2}$. Further, let $\phi_{\varepsilon} * \psi$ be a smooth approximation of $\psi$. Then, by applying Theorem \ref{thm:decay-of-matrix-coeff}, we can deduce the desired result. 
\begin{align*}
\left\vert \left\langle \pi(g_t) f, \frac{\chi_{P_0}}{\mu(P_0)} \right\rangle \right\vert &\leq \int_{SL_2(\R)/\Gamma} \left\vert f(g_t^{-1} g(h\Gamma)) \frac{\chi_{P_0}(g(h\Gamma))}{\mu(P_0)} \right\vert \, d\mu\\ 
&\leq \int_{SL_2(\R)/\Gamma} \left\vert \pi(g_t)f \right\vert \cdot \left\vert \chi - \psi \right\vert \, d\mu + \int_{SL_2(\R)/\Gamma} \left\vert \pi(g_t)f \right\vert \cdot \left\vert \psi - \phi_{\varepsilon}*\psi \right\vert \, d\mu \\
&\,\,\,\,\,\,\,\,\, + \int_{SL_2(\R)/\Gamma} \left\vert \pi(g_t)f \phi_{\varepsilon}*\psi \right\vert \, d\mu\\
& \leq \varepsilon M_f + C b(t) \mathcal{S}_K(f) \mathcal{S}_K(\phi_{\varepsilon}*\psi) 
\end{align*}
where $M_f$ is the maximum of $f$. If $\Gamma$ is tempered, then $b(t)$ is $t e^{-\frac{t}{2}}$ and $C$ is independent of $s$. If $\Gamma$ is non-tempered, $b(t)$ is $e^{-\frac{t}{2}(1-s)}$ and $C$ depends on $s$. It remains to approximate the $K$-Sobolev norm of the function $\phi_{\varepsilon}*\psi$. 
\begin{align*}
\mathcal{S}_K(\phi_{\varepsilon}*\psi)^2 &= \| \phi_{\varepsilon}*\psi \|_2^2 + \| d(\omega) \phi_{\varepsilon}*\psi \|_2^2 \\
& \leq \| \psi \|_2^2 + \frac{1}{\varepsilon^4} \cdot \sup_g \left(d\omega( \phi) \right) \| \psi \|_2^2
\end{align*}
Observe that by an application of Cauchy-Schwarz, and using that $m_P^l(\mathcal{O}_P) = \delta$ we have 
\begin{equation*}
\| \psi \|_2 \leq \| \chi \|_2 + \varepsilon \leq \frac{1}{\delta L} + \varepsilon,
\end{equation*}
so we have, for sufficiently small $\delta$,
\begin{align*}
\left\vert \left\langle \pi(g_t) f, \frac{\chi_{P_0}}{\mu(P_0)} \right\rangle \right\vert & \leq \varepsilon M_f + C b(t) \mathcal{S}_K(f) \left(\| \psi \|_2 \frac{1}{\varepsilon^2} \left( \sqrt{1 + \sup_{g} d\omega(\phi)} \right) \right) \\
& \leq \varepsilon M_f + C b(t) \mathcal{S}_K(f)\left( (\| \chi \|_2 + \varepsilon) \frac{1}{\varepsilon^2} \left( \sqrt{1 + \sup_{g} d\omega(\phi)} \right) \right) \\
& \leq \varepsilon M_f + C b(t) \mathcal{S}_K(f)\left( \left(\frac{1}{\delta L} + \varepsilon \right) \frac{1}{\varepsilon^2} \left( \sqrt{1 + \sup_{g} d\omega(\phi)} \right) \right) \\
& \leq \varepsilon M_f + C b(t) \mathcal{S}_K(f)\left( \left(\frac{2}{\varepsilon^2 \delta L}\right) \left( \sqrt{1 + \sup_{g} d\omega(\phi)} \right) \right) \text{.}
\end{align*}
Hence, combining all of our estimates, we have 
\begin{equation*}
\left\vert \frac{1}{e^t L} \int_{0}^{e^t L} f(h_{s} g_t^{-1} h\Gamma) \, ds \right\vert \leq \varepsilon M_f + C b(t) \mathcal{S}_K(f) \left(\frac{1}{\varepsilon^2 \delta }\right) + D \sqrt{\delta}
\end{equation*}
for some constants $C,D > 0$. To conclude, pick $\varepsilon = e^{-\frac{1}{10}(1-s)}t$ and $\delta = e^{-\frac{1}{5}(1-s)}t$. As $\varepsilon, \delta \to 0$, the uniform convergence of this matrix coefficient to the desired integral gives the result.
\end{proof}

\end{document}